\documentclass[10pt, reqno]{amsart}
\usepackage{amsmath,amssymb,amscd,mathrsfs,amscd}
\usepackage{graphics,verbatim}
\usepackage[all]{xy}

 \newlength{\baseunit}               
 \newcount{\numlines}                
\newtheorem{theorem}{Theorem}
\newtheorem{prop}{Proposition}
\newtheorem{lemma}{Lemma}

\newtheorem*{thm}{Theorem}

\theoremstyle{definition}

\newtheorem{remark}{Remark}

\DeclareMathOperator{\Hom}{Hom} 
 \DeclareMathOperator{\Res}{Res}
\DeclareMathOperator{\End}{End} 
 
\DeclareMathOperator{\I}{I}

\newcommand{\h} {\mathfrak{h}}
\newcommand{\C} {\mathbb{C}}

\newcommand{\Stab}{\mathrm{Stab}}

\newcommand{\af}{\alpha}

\newcommand{\dt}{\delta}
\newcommand{\ep}{\varepsilon}

\newcommand{\ld}{\lambda}

\newcommand{\om}{\omega}

\newcommand{\Dt}{\Delta}

\newcommand{\Q}{{\mathbb{Q}}}
\newcommand{\Z}{{\mathbb{Z}}}

\newcommand{\Sl}{{\mathfrak{sl}}}
\newcommand{\gl}{{\mathfrak{gl}}}
\newcommand{\g}{{\mathfrak{g}}}
\newcommand{\n}{{\mathfrak{n}}}

\newcommand{\s}{{\mathfrak{s}}}

\newcommand{\sgn}{{\operatorname{sgn}}}
\newcommand{\cont}{{\operatorname{cont}}}
\newcommand{\Seq}{{\operatorname{Seq}}}

\newcommand{\E}{{\mathcal{E}}}

\newcommand{\KL}{{\mathcal{KL}}}
\newcommand{\A}{{\mathcal{A}}}
\newcommand{\calC}{{\mathcal{C}}}
\newcommand{\HK}{{\mathcal{HK}}}

\newcommand{\bbE}{{\mathbb{E}}}
\newcommand{\bbK}{{\mathbb{K}}}
\newcommand{\bbI}{{\mathbb{I}}}

\newcommand{\bone}{{\mathbf{1}}}

\newcommand{\ui}{{\underline{i}}}
\newcommand{\uj}{{\underline{j}}}

\newcommand{\uk}{{\underline{k}}}

\renewcommand{\Cap}{{\mbox{$\bigcap$}}}
\renewcommand{\Cup}{{\mbox{$\bigcup$}}}
\newcommand{\Bubble}{{\bigcirc^\bullet}}
\newcommand{\bubble}{{\mbox{\scriptsize$\Bubble$}}}
\newcommand{\yCap}{{\Cap^\bullet}}
\newcommand{\yCup}{{\Cup^\bullet}}
\newcommand{\ycap}{{\mbox{\scriptsize$\yCap$}}}
\newcommand{\ycup}{{\mbox{\scriptsize$\yCup$}}}

\begin{document}
\pagestyle{plain}
\title{The Khovanov-Lauda 2-category and categorifications of a level two quantum sl(n) representation}

\author{David Hill }
\address{Department of Mathematics \\
            University of California, Berkeley \\
            Berkeley, CA 94720-3840}
\email{dhill1@math.berkeley.edu}
\author{Joshua Sussan}
\address{Department of Mathematics \\
            University of California, Berkeley \\
            Berkeley, CA 94720-3840}
\email{sussan@math.berkeley.edu}
\thanks{Research of the authors was partially supported by
NSF EMSW21-RTG grant DMS-0354321}\
\date{\today}

\begin{abstract}
We construct 2-functors from a 2-category categorifying quantum sl(n) to 2-categories categorifying the irreducible representation of highest weight $ 2 \omega_k. $
\end{abstract}

\maketitle



\section{introduction}

Khovanov and Lauda introduced a 2-category whose Grothendieck group is $
\mathcal{U}_q(\mathfrak{sl}_n) $ \cite{kl}. This work generalizes earlier work
by Lauda for the $ \mathcal{U}_q(\mathfrak{sl}_2) $ case \cite{l}. Rouquier has
independently produced a 2-category with similar generators and relations, \cite{r}. There have been
several examples of categorifications of representations of $
\mathcal{U}_q(\mathfrak{sl}_n) $ arising in various contexts. Khovanov and
Lauda conjectured that their 2-category acts on various known categorifications
via a 2-functor.  For example, in their work they construct such a 2-functor to a category
of graded modules over the cohomology of partial flag varieties.  This
2-category categorifies the irreducible representation of $
\mathcal{U}_q(\mathfrak{sl}_n) $ of highest weight $ n \omega_1 $ where $
\omega_1 $ is the first fundamental weight.

In this note we construct this action for the categorification constructed by
Huerfano and Khovanov in \cite{hk}. They categorify the irreducible
representation $ V_{2 \omega_k} $ of highest weight $ 2 \omega_k $, by a
modification of a diagram algebra introduced in \cite{k}. The objects of
2-category $ \mathcal{HK}_{k,n} $ are categories $ \mathcal{C}_{\lambda} $
which are module categories over the modified Khovanov algebra. We explicitly
construct natural transformations between the functors in \cite{hk} and show
that they satisfy the relations in the Khovanov-Lauda 2-category giving the
theorem:

\begin{thm}
There exists a 2-functor $ \Omega_{k,n} \colon \mathcal{KL} \rightarrow
\mathcal{HK}_{k,n}. $
\end{thm}

The Huerfano-Khovanov categorification is based on categories used for the
categorification of $ \mathcal{U}_q(\mathfrak{sl}_2)$-tangle invariants. This
hints that a categorification of $ V_{2 \omega_k} $ may also be obtained on
maximal parabolic subcategories of certain blocks of category $
\mathcal{O}(\mathfrak{gl}_{2k}). $ More specifically, we construct a 2-category
$ \mathcal{P}_{k,n} $ whose objects are full subcategories $ {}_{\mathbb{Z}}
\mathcal{P}^{(k,k)}_{\mu}(\mathfrak{gl}_{2k}) $ of graded
category $
{}_{\mathbb{Z}}\mathcal{O}^{(k,k)}_{\mu}(\mathfrak{gl}_{2k}) $
whose set of objects are those modules which have projective presentations by
projective-injective objects.  The 1-morphisms of $ \mathcal{P}_{k,n} $ are
certain projective functors.  We explicitly construct the 2-morphisms as
natural transformations between the projective functors by the Soergel functor
$ \mathbb{V}.$  We then prove:

\begin{thm}
There is a 2-functor $ \Pi_{k,n} \colon \mathcal{KL} \rightarrow
\mathcal{P}_{k,n}. $
\end{thm}

It should be possible to categorify $ V_{N \omega_k} $ for $ N \geq 1 $ using categories which appear in various knot homologies.  For $
N \geq 2, $ the module categories $ \mathcal{C}_{\lambda} $ in the
Huerfano-Khovanov construction should be replaced by suitable categories of
matrix factorization based on Khovanov-Rozansky link homology. The categories
of matrix factorizations must be generalized from those used in \cite{kr}.
Khovanov and Rozansky suggest that the categories of matrix factorizations
should be taken over tensor products of polynomial rings invariant under the
symmetric group.  These categories were studied in depth by Yonezawa and Wu
\cite{y,w}.  In fact, the isomorphisms of functors categorifying the $
\mathcal{U}_q(\mathfrak{sl}_n) $ relations were defined implicitly in \cite{w}.
To check that there is a a 2-representation of the Khovanov-Lauda 2-category,
these isomorphisms would need to be made more explicit. The category $
\mathcal{O} $ approach should be modified as well.  Now the objects of the
2-category should be subcategories of parabolic subcategories corresponding to the composition $
Nk = k+ \cdots + k $ of blocks of $ \mathcal{O}_{\lambda}(\mathfrak{gl}(Nk)), $
and the stabilizer of the dominant integral weight $ \mu $ is taken to be $ \mathbb{S}_{\lambda_1} \times \cdots \times \mathbb{S}_{\lambda_n} $
where each $ \lambda_i \in \lbrace 0, 1, \ldots, N \rbrace, $ cf. Section~\ref{S:Cat O} below.
Note that a categorification of $ V_{\lambda} $ for arbitrary dominant integral $ \lambda$, hence in particular of $ V_{N\omega_k}$ , is constructed in \cite{bk2} using cyclotomic quotients of Khovanov-Lauda-Rouquier algebras.

While this paper was in preparation, two very relevant papers appeared. 
In \cite{bs}, J. Brundan and C. Stroppel also defined the appropriate natural transformations and checked relations between them to establish a version of the first theorem above, but for Rouquier's 2-category from \cite{r} rather than the Khovanov-Lauda 2-category.
One of the advantages of their result is that they are able to work over an arbitrary field, while we work over a field of characteristic 2. 
It is not immediately clear to us how to use their sign conventions to get an action of the full Khovanov-Lauda 2-category in characteristic zero, because they seem to lead to inconsistencies between propositions \ref{kl3.1kl3.2}, \ref{kl3.3}, \ref{kl3.4b}, and \ref{kl3.10}.
Additionally, Brundan and Stroppel categorify $ V_{2\omega_k} $ using graded category $ \mathcal{O}$. More precisely, they first categorify the classical limit of $ V_{2\omega_k} $ at $q=1$ using a certain parabolic category $ \mathcal{O} $, without mentioning gradings. Then they  establish an equivalence between this category and the (ungraded) diagrammatic category. Finally, they observe that both categories are Koszul (by \cite{bgs} and \cite{bs2}, respectively) so, exploiting unicity of Koszul gradings, their categorification at $q=1$ can be lifted to a categorification of the module $ V_{2\omega_k} $ itself in terms of graded category $ \mathcal{O}$. Our construction on the graded category $ \mathcal{O} $ side is more explicit, relying heavily on the Soergel functor, the Koszul grading that $ \mathcal{O} $ inherits from geometry, and explicit calculations on the cohomology of flag varieties made in \cite{kl}.
In the other relevant paper, M. Mackaay \cite{m} constructs an action of the Khovanov-Lauda 2-category on a category of foams which is the basis of an $ \mathfrak{sl}_3$-knot homology.

{\it Acknowledgements: } The authors would like to thank Mikhail Khovanov and Aaron Lauda for helpful conversations.

\section{The quantum group $ \mathcal{U}_q(\mathfrak{sl}_n) $}

\subsection{Root Data}
Let $\Sl_n=\Sl_n(\C)$ denote the Lie algebra of  traceless $n\times n$-matrices
with standard triangular decomposition $\Sl_n=\n^-\oplus\h\oplus\n^+$. Let
$\Dt\subset\h^*$ be the root system of type $A_{n-1}$ with simple system
$\Pi=\{\af_i|i=1,\ldots,n-1\}$. Let $(\cdot,\cdot)$ denote the symmetric
bilinear form on $\h^*$ satisfying
\[
(\af_i,\af_j)=a_{ij},
\]
where $A=(a_{ij})_{1\leq i,j<n}$ is the Cartan matrix of type $A_{n-1}$:
\begin{eqnarray*}
a_{ij} =\begin{cases} 2 &\mbox{if } j=i, \\
-1 &\mbox{if } |j-i|=1,\\
0 &\mbox{if } |i-j|>1.\end{cases}
\end{eqnarray*}
Let $\Dt^+$ be the set of simple roots relative to $\Pi$. Let
$\om_1,\ldots,\om_{n-1}\in\h^*$ be the elements satisfying
$(\om_i,\af_j)=\dt_{ij}$, and let
$$Q=\bigoplus_{i=1}^{n-1}\Z\af_i,\;\;\;Q^+=\bigoplus_{i=1}^{n-1}\Z_{\geq0}\af_i,
\;\;\;P=\bigoplus_{i=1}^{n-1}\Z\om_i, \;\;\;\mbox{and}\;\;\;
P^+=\bigoplus_{i=1}^{n-1}\Z_{\geq0}\om_i$$ denote the root lattice, positive root
lattice, weight lattice, and dominant weight lattice, respectively.

Set $I=\{1,\ldots,n-1,-1,\ldots,-n+1\}$, $I^+=\{1,\ldots,n-1\}$ and $I^-=-I^+$.
Define $\af_{-i}=-\af_i$, and extend the definition of $a_{ij}$ to all $i,j\in
I$ accordingly. Finally, for $i\in I$, let $\sgn(i)=i/|i|$ be the sign of $i$.

The quantum group $ \mathcal{U}_{q}(\mathfrak{sl}_n) $ is the associative
algebra over $ \mathbb{Q}(q) $ with generators $ E_i, K_i,$ for $i\in I$
satisfying the following conditions:
\begin{enumerate}
\item\label{qrelation1} $ K_i K_{-i} = K_{-i} K_i = 1$, and $ K_i K_j = K_j K_i $ for $i,j\in I$;
\item\label{qrelation2} $ K_i E_j = q^{a_{i,j}} E_j K_i $, $i,j\in I$,
\item\label{qrelation3} $ E_i E_{-j} - E_{-j} E_i = \delta_{i,j} \frac{K_i - K_{-i}}{q-q^{-1}}
$, $i,j\in I^{\pm}$;
\item\label{qrelation4} $ E_i E_j = E_j E_i $, $i,j\in I^{\pm}$, $ |i-j|>1 $;
\item\label{qrelation5} $ E_i^2 E_{j} - (q+q^{-1}) E_i E_{j} E_i + E_{j} E_i^2 = 0 $, $i,j\in I^\pm$, $|i-j|=1$.
\end{enumerate}
We fix a comultiplication $ \Delta \colon \mathcal{U}_q(\mathfrak{sl}_n)
\rightarrow \mathcal{U}_q(\mathfrak{sl}_n) \otimes
\mathcal{U}_q(\mathfrak{sl}_n) $ given as follows for all $i\in I^+$:
\begin{align*}
\Delta(E_i) &= 1 \otimes E_i + E_i \otimes K_i\\
\Delta(E_{-i})&= K_{-i} \otimes E_{-i} + E_{-i} \otimes 1\\
\Delta(K_{\pm i}) &= K_{\pm i} \otimes K_{\pm i}\\
\end{align*}
Via $ \Delta, $ a tensor product of $ \mathcal{U}_q(\mathfrak{sl}_n) $-modules
becomes a $ \mathcal{U}_q(\mathfrak{sl}_n) $-module.

In this paper we are interested in the irreducible $\mathcal{U}_q(\Sl_n)$-modules, $ V_{2 \omega_k} $ with highest weight $ 2 \omega_k. $ Therefore, we will identify the weight lattice $P\cong\Z^{n-1}\subset\Z^n$ as follows: Assume $\ld=\sum_ia_i\om_i$. For each $1\leq i<n$ set
$$\ld_i=\frac{2k-a_1-2a_2- \cdots - (i-1)a_{i-1}+ (n-i)a_i+(n-i-1)a_{i+1}+\cdots+a_{n-1}}{n}$$.

Let $P(2\om_k)$ denote the set of weights of $V_{2\om_k}$. It is well known that under this identification each $\ld\in P(2\om_k)$ satisfies $\ld_i\in\{0,1,2\}$ for all $1\leq i\leq n$ and $\ld_1+\cdots+\ld_n=2k$.
\section{The Khovanov-Lauda 2 category}

Let $\Bbbk$ be a field. The $\Bbbk$-linear 2-category $ \mathcal{KL}$ defined
here was originally constructed in \cite{kl}. The original construction is
defined conveniently in terms of diagrams. We do not present the generators and
relations in terms of diagrams here because it would conflict with the diagrams
used in the construction of the 2-representation in the next section.

Let $ \displaystyle I_\infty=\bigcup_{n\geq0} I^n$, $ \displaystyle I_\infty^+=\bigcup_{n\geq0}(I^+)^n$ where $ I^n $ and $ (I^+)^n $ denote $n$-fold Cartesian products.
Given $\ui=(i_1,i_2,\ldots)\in I_\infty$, let
\[
\cont(\ui)=\sum_{i=1}^{n-1}c_i\af_i,\;\;\;\mbox{where }
c_i=\#\{j|i_j=i\}-\#\{j|i_j=-i\}.
\]
Given $\nu\in Q$, let $\Seq(\nu)=\{\ui\in I_\infty|\cont(\ui)=\nu\}$ and, for
$\nu\in Q^+$, define $\Seq^+(\nu)=\{\ui\in I_\infty^+|\cont(\ui)=\nu\}$.
Finally, define $$\Seq=\bigcup_{\nu\in Q}\Seq(\nu).$$

\subsection{The objects}
The set of objects for this 2-category is the weight lattice, $P$.

\subsection{The 1-morphisms} For each $\ld\in P$, let
$\mathcal{I}_\ld\in\End_{\KL}(\ld)$ be the identity morphism and, for $\ld,\ld'\in P$,
set $\mathcal{I}_\ld\mathcal{I}_\ld'=\dt_{\ld,\ld'}\mathcal{I}_\ld$. For each $i\in I$, we define morphisms
$\E_i\mathcal{I}_\ld\in\Hom_{\KL}(\ld,\ld+\af_i)$. Evidently, we have
$\E_i\mathcal{I}_\ld=\mathcal{I}_{\ld+\af_i}\E_i\mathcal{I}_\ld$. For $\ld,\ld'\in P$, we have
\[
\Hom_{\KL}(\ld,\ld')=\bigoplus_{\substack{\ui\in\Seq\\s\in\Z}}\Bbbk
\mathcal{I}_{\ld'}\E_{\ui}\mathcal{I}_{\ld}\{s\}
\]
where $\E_\ui:=\E_{i_1}\cdots \E_{i_r}$ if $\ui=(i_1,\ldots,i_r)\in I_\infty$,
and $s$ refers to a \emph{grading shift}. Observe that
$\mathcal{I}_{\ld'}\E_{\ui}\mathcal{I}_\ld=0$ unless $\cont(\ui)=\ld'-\ld$, and
$\mathcal{I}_{\ld+\cont(\ui)}\E_\ui\mathcal{I}_\ld=\E_\ui\mathcal{I}_\ld$.

\subsection{The 2-morphisms} The 2-morphisms are generated by
\[
Y_{i;\ld}\in \End_\KL(\E_i \mathcal{I}_\ld),\;\;\;\Psi_{i,j;\ld}\in\Hom_\KL(\E_i\E_j \mathcal{I}_\ld, \E_j\E_i \mathcal{I}_\ld),
\]
\[
\Cup_{i;\ld}\in\Hom_\KL(\mathcal{I}_\ld,\E_{-i} \E_i \mathcal{I}_\ld),
\;\;\;\mbox{and}
\;\;\;\Cap_{i;\ld}\in\Hom_\KL(\E_{-i}\E_i \mathcal{I}_\ld,\mathcal{I}_\ld),
\]
for $i,j\in I^\pm$.
We define $\bone_{i;\ld} \in\End_\KL(\E_i \mathcal{I}_\ld)$ to be the identity transformation.

Let $\ld,\ld'\in P$ and
The degrees of the basic 2
morphisms are given by
\[
\deg Y_{i;\ld}=a_{ii},\;\;\;\deg\Psi_{i,j;\ld}=-a_{ij},\;\;\;
\deg\Cup_{i;\ld}=\deg\Cap_{i;\ld}=1+(\af_i,\ld).
\]

Let $ \ld + \cont(\ui) = \ld+\cont(\uj) = \ld+\cont(\uk)= \ld' $
and
$ \ld' + \cont(\ui') = \ld+\cont(\uj') = \ld''. $
Let $ \Theta_1 \in \Hom_{\KL}(\E_{\ui} \mathcal{I}_{\ld}, \E_{\uj}\mathcal{I}_{\ld}) $ and
$ \Theta_2 \in \Hom_{\KL}(\E_{\ui'} \mathcal{I}_{\ld'}, \E_{\uj'}\mathcal{I}_{\ld'}). $
Then denote the horizontal composition of these 2-morphisms by $ \Theta_2 \Theta_1 $
which is an element of $ \Hom_{\KL}(\E_{\ui'} \mathcal{I}_{\ld'}\E_{\ui} \mathcal{I}_{\ld}, \E_{\uj'} \mathcal{I}_{\ld'}\E_{\uj} \mathcal{I}_{\ld}). $
If $ \Theta_3 \in  \Hom_{\KL}(\E_{\uj} \mathcal{I}_{\ld}, \E_{\uk}\mathcal{I}_{\ld}), $
denote the vertical composition of $ \Theta_3 $ and $ \Theta_1 $ by $ \Theta_3 \circ \Theta_1. $

For convenience of notation, we define the following 2-morphisms. If
$\theta\in\End(\E_\ui \mathcal{I}_{\ld})$ let $ {\theta}^{[j]} = \underbrace{\theta \circ \cdots
\circ\theta}_j. $ For each $i\in I$, define the \emph{bubble}
$$\Bubble^{N}_{i;\ld}=\Cap_{i;\ld}\circ(\bone_{-i;\ld+\af_i}Y_{i;\ld})^{[N]}\circ\Cup_{i;\ld}.$$
Also, define \emph{half bubbles}
$$\yCup_{i;\ld}^N=(\bone_{-i;\ld+\af_i}Y_{i;\ld})^{[N]} \circ
\Cup_{i;\ld}\;\;\;\mbox{and}\;\;\;\yCap_{i;\ld}^N=\Cap_{i;\ld} \circ (Y_{-i;\ld+\af_i} \bone_{i,\ld})^{[N]}.$$

We now define the relations satisfied by these basic 2 morphisms. In what
follows, we omit the argument $\ld$ when the relation is independent of it.

\begin{enumerate}
\item\label{I:sl2relation} $ \mathfrak{sl}_2 $ relations:
\begin{enumerate}
\item\label{I:sl2relation:1} For all $i\in I$,
$$ (\Cap_{-i}  \bone_{i}) \circ (\bone_{i}\Cup_{i}) = \bone_{i} = (\bone_{i}  \Cap_{i}) \circ (\Cup_{-i}
\bone_{i});$$
\item\label{I:sl2relation:2} For all $i\in I^{+}$,
$$ Y_{i} = (\Cap_{-i}\bone_{i}) \circ (\bone_{i}Y_{-i}\bone_{i}) \circ
(\bone_{i}\Cup_{i}) = (\bone_{i}\Cap_{i}) \circ (\bone_{i}Y_{-i}\bone_{i}) \circ
(\Cup_{-i}\bone_{i}). $$
\item\label{I:sl2relation:3} Suppose $i\in I$ and $(-\af_i,\ld)>r+1$, then
$$\Bubble^{r}_{i;\ld} = 0, $$ .
\item\label{I:sl2relation4} Let $i\in I$. If $(\af_i,\ld)\leq-1$,
$$ \Bubble^{-(\af_i,\ld)-1}_{i;\ld}=1.$$
\item\label{I:sl2relation5} Let $i\in I$. If $(\af_i,\ld)\geq 1$, then
\begin{align*}
\bone_{i;\ld-\af_i}  \bone_{-i;\ld} = - {\Psi}_{-i,i;\ld} \circ {\Psi}_{i,-i;\ld}+
\sum_{f=0}^{(\af_i,\ld)-1}
\sum_{g=0}^f \yCup_{-i;\ld}^{[(\af_i,\ld)-f-1]} \circ \Bubble^{[-(\af_i,\ld)-1+g}_{i;\ld}]
 \circ \yCap_{-i;\ld}^{[f-g]}
\end{align*}

\item\label{I:sl2relation6} Let $i\in I^+$. If $(\af_{i},\ld)\leq0$, then
\begin{align*}
(\bone_{i;\ld}  \Cap_{-i;\ld}) \circ (\Psi_{i,i;\ld-\af_i}  \bone_{-i;\ld}) \circ
(\bone_{i;\ld} \Cup_{-i;\ld}) =- \sum_{f=0}^{-(\af_i,\ld)}
{Y_{i;\ld}}^{[-(\af_i,\ld)-f]} \Bubble^{[(\af_i,\ld)-1+f]}_{-i;\ld}.
\end{align*}
If $(\af_i,\ld)\geq -2$, then
\begin{align*}
(\Cap_{i; \ld}  \bone_{i;\ld-\alpha_i}) \circ (\bone_{-i;\ld+\alpha_i}  \Psi_{i,i;\ld-\alpha_i}) \circ (\Cup_{i;\ld}  \bone_{i;\ld-\alpha_i}) =
\sum_{g=0}^{(\af_i,\ld)+2} \Bubble_{i;\ld}^{[-(\af_i,\ld)-1+g]}
{Y_{i;\ld-\af_i}}^{[(\af_i,\ld)-g]}.
\end{align*}
\end{enumerate}

\begin{remark} Note that in \eqref{I:sl2relation5} above the exponent of the bubble may be negative, which is not defined.
To make sense of this, for $ i \in I^+ $, define these symbols
inductively by the formula
\[
\left(\sum_{n\geq
0}\Bubble_{i;\ld}^{(\af_{-i},\ld)-1+n}t^n\right)\left(\sum_{n\geq
0}\Bubble_{-i;\ld}^{(\af_{i},\ld)-1+n}t^n\right)=1
\]
and, $\Bubble_{i;\ld}^{-1}=1$ whenever $(\af_i,\ld)=0$.
\end{remark}
\item\label{I:nilHeckerelations} The nil-Hecke relations:
\begin{enumerate}
\item\label{I:nilHeckerelation1} For each $i\in I^+$, $\Psi_{i,i}^{[2]}=0$.
\item\label{I:nilHeckerelation2}  For $i\in I^+$, $(\Psi_{i,i}  \bone_i)\circ(\bone_i  \Psi_{i,i})\circ(\Psi_{i,i}  \bone_i)
=(\bone_i \Psi_{i,i})\circ(\Psi_{i,i}  \bone_i)\circ (\bone_i\Psi_{i,i}). $
\item\label{I:nilHeckerelation3} For $i\in I^+$, $(\bone_i \bone_i) = (\Psi_{i,i}) \circ (Y_i  \bone_i) - (\bone_i  Y_i) \circ (\Psi_{i,i}) = (Y_i
\bone_i) \circ (\Psi_{i,i}) - (\Psi_{i,i}) \circ (\bone_i  Y_i). $
\item\label{I:nilHeckerelation4} For $j,i\in I^-$,
\begin{align*}
\Psi_{j,i} &= (\Cap_{-j}  \bone_i  \bone_j) \circ (\bone_j  \Cap_{-i}  \bone_{-j}  \bone_i \bone_j) \circ (\bone_j
\bone_i  \Psi_{-j,-i}  \bone_i  \bone_j) \circ(\bone_j  \bone_i  \bone_{-j}  \Cup_i  \bone_j) \circ
(\bone_j  \bone_i  \Cup_j)\\
&=(\bone_i  \bone_j  \Cap_i) \circ (\bone_i  \bone_j  \bone_{-i}  \Cap_j  \bone_i) \circ (\bone_i  \bone_j
\Psi_{-j,-i}  \bone_j  \bone_i) \circ(\bone_i  \Cup_{-j}  \bone_{-i}  \bone_j  \bone_i) \circ (\Cup_{-i}  \bone_j
\bone_i).
\end{align*}
\end{enumerate}
\begin{remark}
For all $i,j\in I^\pm$, set
$ \Psi_{i,-j} = (\bone_{-j}\bone_i\Cap_{-j}) \circ
(\bone_{-j}\Psi_{j,i}\bone_{-j})\circ (\Cup_j\bone_i\bone_{-j}). $
\end{remark}
\item\label{I:Rnurelations} The $ R(\nu) $ relations:
\begin{enumerate}
\item\label{I:nilHeckerelation6} For $i,j\in I^\pm$,
$ (\Psi_{-j,i}) \circ (\Psi_{i,-j}) = \bone_i
 \bone_{-j}. $
\item\label{I:Rnurelations1} For $i,j\in I^+$, $i\neq j$,
$$\Psi_{i,j}\circ\Psi_{j,i}=\begin{cases}\bone_i\bone_j&\mbox{if
}|i-j|>1,\\(i-j)(Y_i\bone_j-\bone_jY_j)&\mbox{if }|i-j|=1.\end{cases}$$
\item\label{I:Rnurelations2} For $i,j\in I^+$, $i\neq j$,
$$(\bone_j  Y_i) \circ (\Psi_{i,j}) = (\Psi_{i,j}) \circ (Y_i
\bone_j),\;\;\;\mbox{and}\;\;\;(Y_j  \bone_i) \circ (\Psi_{i,j}) = (\Psi_{i,j}) \circ
(\bone_i Y_{j}). $$
\item For $i,j,k\in I^+$,
$$(\Psi_{j,k}\bone_i)\circ (\bone_j\Psi_{i,k})\circ(\Psi_{i,j}\bone_k)-(\bone_k\Psi_{i,j})\circ(\Psi_{i,k}\bone_j)\circ(\bone_i \Psi_{j,k})
=\begin{cases}0&i\neq k\mbox{ or }|i-j| = 0,\\ (i-j) \bone_i  \bone_j  \bone_i&i=k\mbox{ and
}|i-j|=1.\end{cases}$$
\end{enumerate}
\end{enumerate}

\section{The Huerfano-Khovanov 2-category}
\label{kh}

\subsection{The Khovanov diagram algebra}
\label{diagramalgebra} Let $ \mathcal{A} = \mathbb{C}[x]/x^2. $ This is a $ \mathbb{Z} $-
graded algebra with multiplication map $ m:\A\otimes\A\rightarrow\A$, such that
$\deg1= -1 $ and $\deg x=1$. There is a comultiplication map $ \Delta \colon
\mathcal{A} \rightarrow \mathcal{A} \otimes \mathcal{A} $ such that $ \Delta(1)
= x \otimes 1 + 1 \otimes x $ and $ \Delta(x) = x \otimes x. $ There is a trace
map $ \text{Tr} \colon \mathcal{A} \rightarrow \mathbb{C} $ such that $
\text{Tr}(x)=1 $ and $ \text{Tr}(1)=0. $ There is also a unit map $ \iota
\colon \mathbb{C} \rightarrow \mathcal{A} $ given by $ \iota(1)=1. $ Also, let
$ \kappa \colon \mathcal{A} \rightarrow \mathcal{A} $ be given by $
\kappa(1)=0, \kappa(x) = 1. $ This algebra gives rise to a two dimensional TQFT
$ \mathfrak{F}$, which is a functor from the category of oriented $1+1$
cobordisms to the category of abelian groups. The functor $ \mathfrak{F} $
sends a disjoint union of $ m $ copies of the circle $ \mathbb{S}^1, $ to $
\mathcal{A}^{\otimes m}. $ For a cobordism $ \mathcal{C}_1, $ from two circles
to one circle,  $ \mathfrak{F}(\mathcal{C}_1) = m. $  For a cobordism $
\mathcal{C}_2 $ from one circle to two circles $ \mathfrak{F}(\mathcal{C}_2) =
\Delta. $ For a cobordism $ \mathcal{C}_3, $ from the empty manifold to $
\mathbb{S}^1, $ $ \mathfrak{F}(\mathcal{C}_3)=\iota. $ For a cobordism $
\mathcal{C}_4 $ from the empty manifold to $ \mathbb{S}^1, $ $
\mathfrak{F}(\mathcal{C}_4)= \text{Tr}. $

\begin{figure}[ht]
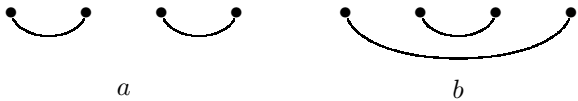

\[
\xy (0,0)*{\bullet};(10,0)*{\bullet}**\crv{(1,-4)&(9,-4)};
(20,0)*{\bullet};(30,0)*{\bullet}**\crv{(21,-4)&(29,-4)};
(15,-10)*{a};
\endxy
\hspace{0.5in}
\xy
(0,0)*{\bullet};(30,0)*{\bullet}**\crv{(2,-8)&(28,-8)};
(10,0)*{\bullet};(20,0)*{\bullet}**\crv{(11,-4)&(19,-4)};
(15,-10)*{b};
\endxy
\]
\caption{Crossingless matches $ a $ and $ b $ for $ r=2$}\label{CM1}
\end{figure}

For any non-negative integer $ r, $  consider $ 2r $ marked points on a line.
Let $ \text{CM}_r $ be the set of non-intersecting curves up to isotopy whose
boundary is the set of the $ 2r $ marked points such that all of the curves lie
on one side of the line. Then there are $ \frac{(2r)!}{r!r!(r+1)} $ elements in
this set.
The set of crossingless matches for $ r=2$ is given in figure ~\ref{CM1}.

\begin{figure}[ht]
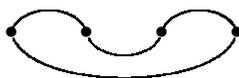

\[
\xy
(0,0)*{};(10,0)*{}**\crv{(1,4)&(9,4)};
(20,0)*{};(30,0)*{}**\crv{(21,4)&(29,4)};
(0,0)*{\bullet};(30,0)*{\bullet}**\crv{(2,-8)&(28,-8)};
(10,0)*{\bullet};(20,0)*{\bullet}**\crv{(11,-4)&(19,-4)};
\endxy
\]
\caption{Concatenation $ (Ra)b$}\label{CM2}
\end{figure}

Let $ a, b \in \text{CM}_r. $  Then $ (Rb)a $ is a collection of
circles obtained by concatenating $ a \in \text{CM}_r $ with the reflection $Rb$
of $ b \in \text{CM}_r $ in the line. Then applying the two dimensional TQFT $
\mathfrak{F}, $ one associates the graded vector space $ {}_bH^{r}_{a} $ to
this collection of circles.  Taking direct sums over all crossingless matches
gives a graded vector space
$$H^{r} = \bigoplus_{a,b}\, {}_{b} H^{r}_{a} \lbrace r \rbrace $$
where the degree $ i $ component of $  {}_{b} H^{r}_{a} \lbrace r \rbrace $ is the degree $ i-r $ component of
$ {}_{b} H^{r}_{a}. $
This graded vector space obtains the structure of an associative algebra via $
\mathfrak{F} $, cf. \cite{k}.

\begin{figure}[ht]
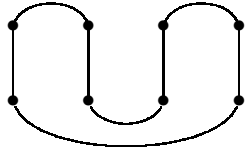

\[
\xy
(0,10)*{};(10,10)*{}**\crv{(1,14)&(9,14)};
(20,10)*{};(30,10)*{}**\crv{(21,14)&(29,14)};
(0,10)*{\bullet};(30,10)*{\bullet};
(10,10)*{\bullet};(20,10)*{\bullet};
(0,0)*{\bullet};(30,0)*{\bullet}**\crv{(2,-8)&(28,-8)};
(10,0)*{\bullet};(20,0)*{\bullet}**\crv{(11,-4)&(19,-4)};
(0,0)*{};(0,10)*{}**\dir{-};
(10,0)*{};(10,10)*{}**\dir{-};
(20,0)*{};(20,10)*{}**\dir{-};
(30,0)*{};(30,10)*{}**\dir{-};
\endxy
\]
\caption{Concatenation $ {}_a T_{b}$}\label{RaTb}
\end{figure}

Let $ T $ be a tangle from $ 2r $ points to $ 2s $ points.  Let $ a $ be a crossingless match for $ 2s $ points and $ b $ a crossingless match for $ 2s $ points.
Then let $ {}_a T_b $ be the concatenation $ Ra \circ T \circ b $ and $ {}_a \mathfrak{F}(T)_b = \mathfrak{F}({}_a T_b).$  See figure ~\ref{RaTb} for an example when $ T $ is the identity tangle.

To any tangle diagram $ T $ from $ 2r $ points to $ 2s $ points, there is a
$(H^s, H^r)$-bimodule
$$ \mathfrak{F}(T)= \bigoplus_{\substack{a \in \text{CM}_r \\
b \in \text{CM}_s}}\, \mathfrak{F}({}_{a} T_{b}) \lbrace r \rbrace . $$
To any cobordism $ C $ between
tangle $ T_1 $ and $ T_2$, there is a bimodule map $ \mathfrak{F}(C) \colon
\mathfrak{F}(T_1) \rightarrow \mathfrak{F}(T_2), $
of degree $ - \chi(C)-r-s $ where $ \chi(C) $ is the Euler characteristic of $ C $
cf. proposition 5 of  \cite{k}.

\begin{figure}[ht]
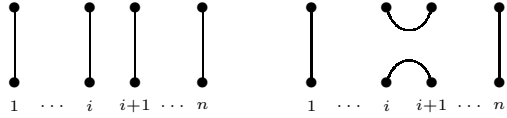

\[
\xy (0,0)*{\bullet};(0,10)*{}**\dir{-};(0,10)*{\bullet};
(10,0)*{\bullet};(10,10)*{}**\dir{-};(10,10)*{\bullet};(16,0)*{\bullet};(16,10)*{}**\dir{-};(16,10)*{\bullet};
(25,0)*{\bullet};(25,10)*{}**\dir{-};(25,10)*{\bullet}; (0,-3)*{\mbox{\tiny
1}};(5,-3)*{\mbox{\tiny $\ldots$}};(10,-3)*{\mbox{\tiny
$i$}};(16,-3)*{\mbox{\tiny
$i$+$1$}};(21,-3)*{\mbox{\tiny$\ldots$}};(25,-3)*{\mbox{\tiny$n$}};
\endxy
\hspace{.5in} \xy (0,0)*{\bullet};(0,10)*{}**\dir{-};(0,10)*{\bullet};
(10,0)*{\bullet};(16,0)*{\bullet}**\crv{(11,4)&(15,4)};
(10,10)*{};(16,10)*{}**\crv{(11,6)&(15,6)};(10,10)*{\bullet};(16,10)*{\bullet};
(25,0)*{\bullet};(25,10)*{}**\dir{-};(25,10)*{\bullet}; (0,-3)*{\mbox{\tiny
1}};(5,-3)*{\mbox{\tiny $\ldots$}};(10,-3)*{\mbox{\tiny
$i$}};(16,-3)*{\mbox{\tiny
$i$+$1$}};(21,-3)*{\mbox{\tiny$\ldots$}};(25,-3)*{\mbox{\tiny$n$}};
\endxy
\]
\caption{I and $U_i$}\label{d1}
\end{figure}

Consider the tangles $ \I$ and $ U_i $ in figure ~\ref{d1}.   Then there are saddle cobordisms
$ S_i \colon U_i \rightarrow \I $ and
$ S^i \colon \I \rightarrow U_i. $

\begin{figure}[ht]
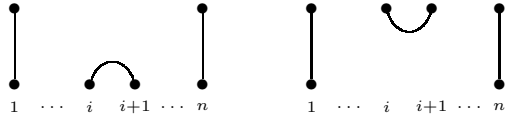

\[
\xy (0,0)*{\bullet};(0,10)*{}**\dir{-};(0,10)*{\bullet};
(10,0)*{\bullet};(16,0)*{\bullet}**\crv{(11,4)&(15,4)};
(25,0)*{\bullet};(25,10)*{}**\dir{-};(25,10)*{\bullet}; (0,-3)*{\mbox{\tiny
1}};(5,-3)*{\mbox{\tiny $\ldots$}};(10,-3)*{\mbox{\tiny
$i$}};(16,-3)*{\mbox{\tiny
$i$+$1$}};(21,-3)*{\mbox{\tiny$\ldots$}};(25,-3)*{\mbox{\tiny$n$}};
\endxy
\hspace{0.5in} \xy (0,0)*{\bullet};(0,10)*{}**\dir{-};(0,10)*{\bullet};
(10,10)*{};(16,10)*{}**\crv{(11,6)&(15,6)};(10,10)*{\bullet};(16,10)*{\bullet};
(25,0)*{\bullet};(25,10)*{}**\dir{-};(25,10)*{\bullet}; (0,-3)*{\mbox{\tiny
1}};(5,-3)*{\mbox{\tiny $\ldots$}};(10,-3)*{\mbox{\tiny
$i$}};(16,-3)*{\mbox{\tiny
$i$+$1$}};(21,-3)*{\mbox{\tiny$\ldots$}};(25,-3)*{\mbox{\tiny$n$}};
\endxy
\]
\caption{$T_i$ and $T^i$}\label{d2}
\end{figure}

\begin{lemma}
\label{mumap}
Let $ T_i $ and $ T^i $ be the tangles in Figure~\ref{d2}.
\begin{enumerate}
\item There exists an $ (H^{n-1},H^n)$-bimodule homomorphism $\mu_i\colon\mathfrak{F}(T_i)\rightarrow \mathfrak{F}(T_{i+1})$ of degree one.
\item There exists an $ (H^{n},H^{n-1})$-bimodule homomorphism $\mu^i\colon\mathfrak{F}(T^i)\rightarrow\mathfrak{F}(T^{i+1})$ of degree one.
\end{enumerate}
\end{lemma}

\begin{proof}
There is a degree zero isomorphism of bimodules
$ \mathfrak{F}(T_i) \cong \mathfrak{F}(T_i) \otimes_{H^n}  \mathfrak{F}(\I) . $
Then by \cite{k} there is a bimodule map of degree one
$$ 1 \otimes \mathfrak{F}(S^{i+1})
\colon \mathfrak{F}(T_i) \otimes_{H^n} \mathfrak{F}(\I)\rightarrow
\mathfrak{F}(T_i) \otimes_{H^n} \mathfrak{F}(U_{i+1}) $$
where $ 1 $ denotes the identity map.
Finally note $
\mathfrak{F}(T_i) \otimes_{H^n}  \mathfrak{F}(U_{i+1}) \cong
\mathfrak{F}(T_{i+1}). $ Then $ \mu_i $ is the composition of these maps.

The construction of $ \mu^i $ is similar.
\end{proof}

\begin{lemma}
\label{mulemma1}
Let $ a \in \text{CM}_n $ and $ b \in \text{CM}_{n-1} $ be two crossingless matches.
Let $ T^i $ be the tangle on the right side of the Figure~\ref{d2}.  Let $ U_i $ be the tangle in Figure~\ref{d1}. Consider the
homomorphism induced by the cobordism $ S^i, $ $ \mathfrak{F}(T^i) \rightarrow \mathfrak{F}(U_i) \otimes_{H^n}
\mathfrak{F}(T^i) \cong \A \otimes_{\mathbb{C}}
\mathfrak{F}(T^i) $. Let $ \alpha \otimes \beta \in \mathfrak{F}({}_a T^i{}_b) $
where $ \alpha \in \mathcal{A} $ corresponds to the circle passing through the
point $ i $ on the top line and $ \beta \in \mathcal{A}^{\otimes p} $ corresponds to the
remaining circles. Then $ \alpha \otimes \beta \mapsto \Delta(\alpha) \otimes
\beta. $
\end{lemma}

\begin{proof}
The map is induced by the cobordism $ S^i. $ On the set of circles, this
cobordism is a union of identity cobordisms and a cobordism $ \mathcal{C}_2. $
The result now follows upon applying $ \mathfrak{F}. $

\end{proof}

\begin{lemma}
\label{mulemma2} Let $ \I $ be the identity tangle from $ 2r $ points to $ 2r $
points, $ T_i $ a tangle from $ 2(r+1) $ points to $ 2r $ points and $ T^i $ a
tangle from $ 2r $ points to $ 2(r+1) $ points. Let $ a $ and $ b $ be cup
diagrams for $ 2r $ points. Consider the map
$$ \A \otimes_{\mathbb{C}} \mathfrak{F}(\I)
\rightarrow \mathfrak{F}(T_i) \otimes_{H^{r+1}} \mathfrak{F}(T^i) \rightarrow
\mathfrak{F}(T_{i+1})\otimes_{H^{r+1}} \mathfrak{F}(T^i) \rightarrow \mathfrak{F}(\I) $$
where the first and last maps are isomorphisms and the middle map is $ \mu_i
\otimes 1. $ Let $ \beta \in \mathcal{A} $ correspond to the circle passing
through point $ i $ of $ {}_{a}\I{}_{b}, $ $ \gamma \in \mathcal{A}^{\otimes r} $
correspond to the remaining circles and $ \alpha \in \mathcal{A}. $ Then the
map above sends $ \alpha \otimes \beta \otimes \gamma \mapsto (\alpha \beta)
\otimes \gamma. $
\end{lemma}

\begin{proof}
The map is induced by a cobordism $ S^{i+1}. $ On the set of circles, this
cobordism is union of identity cobordisms and a cobordism $ \mathcal{C}_1. $
The result now follows upon applying $ \mathfrak{F}. $
\end{proof}

\subsection{The Huerfano-Khovanov categorification}
\label{HKcategory} Let  $ \lambda \in P(2\om_k)$.
Recall that
$\af_{-i}=-\af_i$. Hence, for $i\in I$, we have
$$ \lambda + \alpha_i =
(\lambda_1, \ldots, \lambda_i +\sgn(i), \lambda_{i+1}-\sgn(i),\ldots,
\lambda_n).$$
Label $ n $ collinear points by the integers $ \lambda_i. $  Those
points labeled by $ 0 $ or $ 2 $ will never be the boundaries of arcs but will
rather just serve as place holders. Then define the algebra $ H_{\lambda} =
H^{\gamma(\lambda)} $ where $ \gamma(\lambda) = \frac{1}{2}|\lbrace \lambda_i |
\lambda_i =1 \rbrace|. $  Let $ e_{\ld} $ be the identity element.

Let $i\in I^+$. We define five special tangles $D_{\lambda,i}, D^{\lambda, i},
T_{\lambda,i}, T^{\lambda, i}, \text{I}_{\ld} $ in figures ~\ref{d3}, ~\ref{d4}, ~\ref{d5}. If
a point is labeled by zero or two, it will not be part of the boundary of any
curve.  Away from points $ i, i+1 $ the tangle is the identity.

\begin{figure}[ht]
\[
\xy (0,0)*{\bullet};(0,10)*{}**\dir{-};(0,10)*{\bullet};
(10,0)*{};(16,10)*{}**\dir{-};(16,10)*{\bullet};(10,10)*{\bullet};(16,0)*{\bullet};(10,0)*{\bullet};
(25,0)*{\bullet};(25,10)*{}**\dir{-};(25,10)*{\bullet}; (0,-3)*{\mbox{\tiny
$\ld_1$}};(5,-3)*{\mbox{\tiny $\ldots$}};(10,-3)*{\mbox{\tiny
$\ld_i$}};(16,-3)*{\mbox{\tiny
$\ld_{i+1}$}};(21,-3)*{\mbox{\tiny$\ldots$}};(26,-3)*{\mbox{\tiny$\ld_n$}};
(0,13)*{\mbox{\tiny
$\ld_1$}};(4,13)*{\mbox{\tiny $\ldots$}};(9,13)*{\mbox{\tiny
$\ld_i$+1}};(16,13)*{\mbox{\tiny
$\ld_{i+1}$-1}};(22,13)*{\mbox{\tiny$\ldots$}};(26,13)*{\mbox{\tiny$\ld_n$}};
\endxy
\hspace{.5in} \xy (0,0)*{\bullet};(0,10)*{}**\dir{-};(0,10)*{\bullet};
(10,0)*{\bullet};(16,10)*{\bullet};(10,10)*{\bullet};(16,0)*{};(10,10)*{}**\dir{-};(16,0)*{\bullet};
(25,0)*{\bullet};(25,10)*{}**\dir{-};(25,10)*{\bullet}; (0,-3)*{\mbox{\tiny
$\ld_1$}};(5,-3)*{\mbox{\tiny $\ldots$}};(10,-3)*{\mbox{\tiny
$\ld_i$}};(16,-3)*{\mbox{\tiny
$\ld_{i+1}$}};(21,-3)*{\mbox{\tiny$\ldots$}};(26,-3)*{\mbox{\tiny$\ld_n$}};
(0,13)*{\mbox{\tiny
$\ld_1$}};(4,13)*{\mbox{\tiny $\ldots$}};(9,13)*{\mbox{\tiny
$\ld_i$+1}};(16,13)*{\mbox{\tiny
$\ld_{i+1}$-1}};(22,13)*{\mbox{\tiny$\ldots$}};(26,13)*{\mbox{\tiny$\ld_n$}};
\endxy
\]
\caption{$D_{\ld,i}$ and $D^{\ld,i}$}\label{d3}
\end{figure}

\begin{figure}[ht]
\[
\xy (0,0)*{\bullet};(0,10)*{}**\dir{-};(0,10)*{\bullet};(16,10)*{\bullet};(10,10)*{\bullet};
(10,0)*{\bullet};(16,0)*{\bullet}**\crv{(11,4)&(15,4)};
(25,0)*{\bullet};(25,10)*{}**\dir{-};(25,10)*{\bullet}; (0,-3)*{\mbox{\tiny
$\ld_1$}};(5,-3)*{\mbox{\tiny $\ldots$}};(10,-3)*{\mbox{\tiny
$\ld_i$}};(16,-3)*{\mbox{\tiny
$\ld_{i+1}$}};(21,-3)*{\mbox{\tiny$\ldots$}};(25,-3)*{\mbox{\tiny$\ld_n$}};
(0,13)*{\mbox{\tiny
$\ld_1$}};(4,13)*{\mbox{\tiny $\ldots$}};(9,13)*{\mbox{\tiny
$\ld_i$+1}};(16,13)*{\mbox{\tiny
$\ld_{i+1}$-1}};(22,13)*{\mbox{\tiny$\ldots$}};(26,13)*{\mbox{\tiny$\ld_n$}};
\endxy
\hspace{0.5in} \xy (0,0)*{\bullet};(0,10)*{}**\dir{-};(0,10)*{\bullet};(16,0)*{\bullet};(10,0)*{\bullet};
(10,10)*{};(16,10)*{}**\crv{(11,6)&(15,6)};(10,10)*{\bullet};(16,10)*{\bullet};
(25,0)*{};(25,10)*{}**\dir{-};(25,10)*{\bullet}; (0,-3)*{\mbox{\tiny
$\ld_1$}};(5,-3)*{\mbox{\tiny $\ldots$}};(10,-3)*{\mbox{\tiny
$\ld_i$}};(16,-3)*{\mbox{\tiny
$\ld_{i+1}$}};(21,-3)*{\mbox{\tiny$\ldots$}};(25,-3)*{\mbox{\tiny$\ld_n$}};
(0,13)*{\mbox{\tiny
$\ld_1$}};(4,13)*{\mbox{\tiny $\ldots$}};(9,13)*{\mbox{\tiny
$\ld_i$+1}};(16,13)*{\mbox{\tiny
$\ld_{i+1}$-1}};(22,13)*{\mbox{\tiny$\ldots$}};(26,13)*{\mbox{\tiny$\ld_n$}};
\endxy
\]
\caption{$T_{\ld,i}$ and $T^{\ld,i}$}\label{d4}
\end{figure}

\begin{figure}[ht]
\[
\xy (0,0)*{\bullet};(0,10)*{}**\dir{-};(0,10)*{\bullet};
(10,0)*{\bullet};(10,10)*{}**\dir{-};(10,10)*{\bullet};(16,0)*{\bullet};(16,10)*{}**\dir{-};(16,10)*{\bullet};
(25,0)*{\bullet};(25,10)*{}**\dir{-};(25,10)*{\bullet}; (0,-3)*{\mbox{\tiny
$\ld_1$}};(5,-3)*{\mbox{\tiny $\ldots$}};(10,-3)*{\mbox{\tiny
$\ld_i$}};(16,-3)*{\mbox{\tiny
$\ld_{i+1}$}};(21,-3)*{\mbox{\tiny$\ldots$}};(25,-3)*{\mbox{\tiny$\ld_n$}};
(0,13)*{\mbox{\tiny
$\ld_1$}};(5,13)*{\mbox{\tiny $\ldots$}};(10,13)*{\mbox{\tiny
$\ld_i$}};(16,13)*{\mbox{\tiny
$\ld_{i+1}$}};(21,13)*{\mbox{\tiny$\ldots$}};(25,13)*{\mbox{\tiny$\ld_n$}};
\endxy
\]
\caption{Identity tangle $\text{I}_\ld$}\label{d5}
\end{figure}

\begin{figure}[ht]
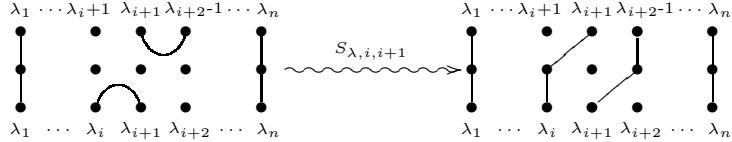

\[
\xy (0,0)*{ \xy (0,0)*{\bullet};(0,10)*{}**\dir{-};(0,10)*{\bullet};(10,10)*{\bullet};
(10,0)*{\bullet};(16,0)*{\bullet}**\crv{(11,4)&(15,4)};(22,0)*{\bullet};
(16,10)*{};(22,10)*{}**\crv{(17,6)&(21,6)};(16,10)*{\bullet};(22,10)*{\bullet};
(32,0)*{\bullet};(32,10)*{}**\dir{-};(32,10)*{\bullet};
(0,5)*{\bullet};(10,5)*{\bullet};(16,5)*{\bullet};(22,5)*{\bullet};(32,5)*{\bullet};
(0,-3)*{\mbox{\tiny $\ld_1$}};(5,-3)*{\mbox{\tiny
$\ldots$}};(10,-3)*{\mbox{\tiny $\ld_i$}};(16,-3)*{\mbox{\tiny
$\ld_{i+1}$}};(22.5,-3)*{\mbox{\tiny$\ld_{i+2}$}};(28,-3)*{\mbox{\tiny$\ldots$}};(33,-3)*{\mbox{\tiny$\ld_n$}};
(0,13)*{\mbox{\tiny $\ld_1$}};(4,13)*{\mbox{\tiny
$\ldots$}};(9,13)*{\mbox{\tiny $\ld_i$+1}};(16,13)*{\mbox{\tiny
$\ld_{i+1}$}};(23,13)*{\mbox{\tiny$\ld_{i+2}$-1}};(29,13)*{\mbox{\tiny$\ldots$}};(33,13)*{\mbox{\tiny$\ld_n$}};
\endxy
}="x"; (60,0)*{ \xy (0,0)*{\bullet};(0,10)*{}**\dir{-};(0,10)*{\bullet};
(10,0)*{\bullet};(10,5)*{}**\dir{-};(16,10)*{}**\dir{-};(22,0)*{\bullet};(16,0)*{\bullet};
(10,10)*{\bullet};
(16,0)*{\bullet};(22,5)*{}**\dir{-};(22,10)*{\bullet}**\dir{-};(16,10)*{\bullet};(22,10)*{\bullet};
(32,0)*{\bullet};(32,10)*{}**\dir{-};(32,10)*{\bullet};
(0,5)*{\bullet};(10,5)*{\bullet};(16,5)*{\bullet};(22,5)*{\bullet};(32,5)*{\bullet};
(0,-3)*{\mbox{\tiny $\ld_1$}};(5,-3)*{\mbox{\tiny
$\ldots$}};(10,-3)*{\mbox{\tiny $\ld_i$}};(16,-3)*{\mbox{\tiny
$\ld_{i+1}$}};(22.5,-3)*{\mbox{\tiny$\ld_{i+2}$}};(28,-3)*{\mbox{\tiny$\ldots$}};(33,-3)*{\mbox{\tiny$\ld_n$}};
(0,13)*{\mbox{\tiny $\ld_1$}};(4,13)*{\mbox{\tiny
$\ldots$}};(9,13)*{\mbox{\tiny $\ld_i$+1}};(16,13)*{\mbox{\tiny
$\ld_{i+1}$}};(23,13)*{\mbox{\tiny$\ld_{i+2}$-1}};(29,13)*{\mbox{\tiny$\ldots$}};(33,13)*{\mbox{\tiny$\ld_n$}};
\endxy
}="y"; {\ar@{~>}^{\mbox{\tiny$S_{\ld,i,i+1}$}}"x";"y"};
\endxy
\]
\caption{Cobordism $S_{\ld,i,i+1}$}\label{d6}
\end{figure}

The cobordisms $ S_{\lambda,i} \colon T^{\lambda+\alpha_i,i} \circ
T_{\lambda,i} \rightarrow \I_{\lambda} $ and $ S_{\lambda,i,j} \colon
T^{\lambda+\alpha_{i},j} \circ T_{\lambda,i} \rightarrow D_{\lambda+\alpha_j,i}
\circ D_{\lambda,j}  $ are saddle cobordisms for $ j=i \pm 1. $ Similarly, the
cobordisms $ S^{\lambda,i}, S^{\lambda,i,j} $ are saddle cobordism in the
opposite direction. For example, the cobordism $ S_{\lambda,i,i+1} $ is given
in figure ~\ref{d6}.

Let $ \mathcal{C}_{\lambda} $ be the category of finitely generated, graded $
H_{\lambda}$-modules, and let $\bbI_\ld:\calC_\ld\rightarrow\calC_\ld$ be the
identity functor. For $\ld,\ld'\in P(2\om_k)$, set
$\bbI_{\ld'}\bbI_\ld=\dt_{\ld,\ld'}\bbI_\ld$.

Let $i\in I^+$. To make future definitions more homogeneous, define
$D_{\lambda,-i}$, $D^{\lambda, -i}$,
$T_{\lambda,-i}$, $T^{\lambda, -i}$ as in figures ~\ref{d3a} and ~\ref{d4a}. Also, in what follows, interpret the pair $(\ld_{-i},\ld_{-i+1})$ as $(\ld_{i+1},\ld_i)$ and recall that $\af_{-i}=-\af_i$.

\begin{figure}[ht]
\[
\xy (0,0)*{\bullet};(0,10)*{}**\dir{-};(0,10)*{\bullet};
(10,0)*{};(16,10)*{}**\dir{-};(16,10)*{\bullet};(10,10)*{\bullet};(16,0)*{\bullet};(10,0)*{\bullet};
(25,0)*{\bullet};(25,10)*{}**\dir{-};(25,10)*{\bullet}; (0,-3)*{\mbox{\tiny
$\ld_1$}};(5,-3)*{\mbox{\tiny $\ldots$}};(10,-3)*{\mbox{\tiny
$\ld_i$}};(16,-3)*{\mbox{\tiny
$\ld_{i+1}$}};(21,-3)*{\mbox{\tiny$\ldots$}};(26,-3)*{\mbox{\tiny$\ld_n$}};
(0,13)*{\mbox{\tiny
$\ld_1$}};(4,13)*{\mbox{\tiny $\ldots$}};(9,13)*{\mbox{\tiny
$\ld_i$-1}};(16,13)*{\mbox{\tiny
$\ld_{i+1}$+1}};(22,13)*{\mbox{\tiny$\ldots$}};(26,13)*{\mbox{\tiny$\ld_n$}};
\endxy
\hspace{.5in} \xy (0,0)*{\bullet};(0,10)*{}**\dir{-};(0,10)*{\bullet};
(10,0)*{\bullet};(16,10)*{\bullet};(10,10)*{\bullet};(16,0)*{};(10,10)*{}**\dir{-};(16,0)*{\bullet};
(25,0)*{\bullet};(25,10)*{}**\dir{-};(25,10)*{\bullet}; (0,-3)*{\mbox{\tiny
$\ld_1$}};(5,-3)*{\mbox{\tiny $\ldots$}};(10,-3)*{\mbox{\tiny
$\ld_i$}};(16,-3)*{\mbox{\tiny
$\ld_{i+1}$}};(21,-3)*{\mbox{\tiny$\ldots$}};(26,-3)*{\mbox{\tiny$\ld_n$}};
(0,13)*{\mbox{\tiny
$\ld_1$}};(4,13)*{\mbox{\tiny $\ldots$}};(9,13)*{\mbox{\tiny
$\ld_i$-1}};(16,13)*{\mbox{\tiny
$\ld_{i+1}$+1}};(22,13)*{\mbox{\tiny$\ldots$}};(26,13)*{\mbox{\tiny$\ld_n$}};
\endxy
\]
\caption{$D^{\ld,-i}$ and $D_{\ld,-i}$}\label{d3a}
\end{figure}

\begin{figure}[ht]
\[
\xy (0,0)*{\bullet};(0,10)*{}**\dir{-};(0,10)*{\bullet};(16,10)*{\bullet};(10,10)*{\bullet};
(10,0)*{\bullet};(16,0)*{\bullet}**\crv{(11,4)&(15,4)};
(25,0)*{\bullet};(25,10)*{}**\dir{-};(25,10)*{\bullet}; (0,-3)*{\mbox{\tiny
$\ld_1$}};(5,-3)*{\mbox{\tiny $\ldots$}};(10,-3)*{\mbox{\tiny
$\ld_i$}};(16,-3)*{\mbox{\tiny
$\ld_{i+1}$}};(21,-3)*{\mbox{\tiny$\ldots$}};(25,-3)*{\mbox{\tiny$\ld_n$}};
(0,13)*{\mbox{\tiny
$\ld_1$}};(4,13)*{\mbox{\tiny $\ldots$}};(9,13)*{\mbox{\tiny
$\ld_i$-1}};(16,13)*{\mbox{\tiny
$\ld_{i+1}$+1}};(22,13)*{\mbox{\tiny$\ldots$}};(26,13)*{\mbox{\tiny$\ld_n$}};
\endxy
\hspace{0.5in} \xy (0,0)*{\bullet};(0,10)*{}**\dir{-};(0,10)*{\bullet};(16,0)*{\bullet};(10,0)*{\bullet};
(10,10)*{};(16,10)*{}**\crv{(11,6)&(15,6)};(10,10)*{\bullet};(16,10)*{\bullet};
(25,0)*{};(25,10)*{}**\dir{-};(25,10)*{\bullet}; (0,-3)*{\mbox{\tiny
$\ld_1$}};(5,-3)*{\mbox{\tiny $\ldots$}};(10,-3)*{\mbox{\tiny
$\ld_i$}};(16,-3)*{\mbox{\tiny
$\ld_{i+1}$}};(21,-3)*{\mbox{\tiny$\ldots$}};(25,-3)*{\mbox{\tiny$\ld_n$}};
(0,13)*{\mbox{\tiny
$\ld_1$}};(4,13)*{\mbox{\tiny $\ldots$}};(9,13)*{\mbox{\tiny
$\ld_i$-1}};(16,13)*{\mbox{\tiny
$\ld_{i+1}$+1}};(22,13)*{\mbox{\tiny$\ldots$}};(26,13)*{\mbox{\tiny$\ld_n$}};
\endxy
\]
\caption{$T_{\ld,-i}$ and $T^{\ld,-i}$}\label{d4a}
\end{figure}

Let $i\in I$. Let $\bbI_\ld:\calC_\ld\to\calC_\ld$ denote the identity functor which is tensoring with the $ (H_{\ld}, H_{\ld})- $ bimodule $ H_{\ld}$. Let $\mathbb{E}_i\bbI_\ld \colon \mathcal{C}_{\lambda}
\rightarrow \mathcal{C}_{\lambda + \alpha_i} $ be the functor of tensoring with
a bimodule defined as follows:
\[
\bbE_i\bbI_\ld=\begin{cases}
\mathfrak{F}(D_{\lambda, i}) &\text{if } (\lambda_i, \lambda_{i+1})=(1,2)\\
\mathfrak{F}(D^{\lambda, i}) &\text{if } (\lambda_i, \lambda_{i+1})=(0,1)\\
\mathfrak{F}(T_{\lambda, i}) &\text{if } (\lambda_i, \lambda_{i+1})=(1,1)\\
\mathfrak{F}(T^{\lambda,i}) &\text{if } (\lambda_i, \lambda_{i+1})=(0,2)\\
0&\mbox{otherwise.}
\end{cases}
\]
Evidently, $\bbE_i\bbI_\ld=\bbI_{\ld+\af_i}\bbE_i\bbI_\ld$
for all $i\in I$, and $\bbI_\ld=\mathfrak{F}(\I_\ld)$.


For $i\in I$, let $ \mathbb{K}_i\bbI_\ld \colon \mathcal{C}_{\lambda}
\rightarrow \mathcal{C}_{\lambda}$ be the grading shift functor
$\bbK_i\bbI_\ld=\bbI_\ld\{(\af_i,\ld)\}.$ Finally, set $\calC=\bigoplus_{\ld\in
P(2\om_k)}\calC_\ld$, $\bbE_i=\bigoplus_{\ld\in P(2\om_k)}\bbE_i\bbI_\ld$,
$\bbK_i=\bigoplus_{\ld\in P(2\om_k)}\bbK_i\bbI_\ld$, and
$\bbI=\bigoplus_{\ld\in P(2\om_k)}\bbI_\ld$.

Propositions 2 and 3 of \cite{hk} are that these functors satisfy quantum
$\mathfrak{sl}_n$ relations:
\begin{prop}\cite[Proposition 2,3]{hk} We have
\begin{enumerate}
\item $ \mathbb{K}_i \mathbb{K}_{-i}\bbI_{\ld} \cong \bbI_{\ld}  \cong \mathbb{K}_{-i} \mathbb{K}_i\bbI_{\ld} $, and
$\mathbb{K}_i \mathbb{K}_j \bbI_{\ld} \cong \mathbb{K}_j \mathbb{K}_i\bbI_{\ld}$ for $i,j\in I$;
\item $ \mathbb{K}_i \mathbb{E}_j \bbI_{\ld} \cong \mathbb{E}_j \mathbb{K}_i\bbI_{\ld} \{a_{ij}\}$,
for $i,j\in I$;
\item $ \mathbb{E}_i \mathbb{E}_{-j}\bbI_{\ld} \cong \mathbb{E}_{-j} \mathbb{E}_i\bbI_{\ld} $ if $i,j\in I^+$, $ i \neq j; $
\item $ \mathbb{E}_i \mathbb{E}_j\bbI_{\ld} \cong \mathbb{E}_j \mathbb{E}_i\bbI_{\ld} $ if $i,j\in I^\pm$, $ |i-j|>1; $
\item $ \mathbb{E}_i \mathbb{E}_i \mathbb{E}_j\bbI_{\ld} \oplus \mathbb{E}_j \mathbb{E}_i\mathbb{E}_i\bbI_{\ld} \cong
\mathbb{E}_i \mathbb{E}_j \mathbb{E}_i\bbI_{\ld} \{1\} \oplus \mathbb{E}_i \mathbb{E}_j
\mathbb{E}_i\bbI_{\ld} \{-1\}$ if $i,j\in I^{\pm}$,  $|i-j|=1;$
\item For $i\in I$,
    $$\mathbb{E}_i \mathbb{E}_{-i}\bbI_\ld
    \cong\begin{cases}\mathbb{E}_{-i} \mathbb{E}_i\bbI_\ld \oplus \bbI_\ld \{1\} \oplus \bbI_\ld
    \{-1\}&\mbox{if }i\in I^+\mbox{ and }(\lambda_i, \lambda_{i+1})=(2,0),\\
    \mathbb{E}_{-i} \mathbb{E}_i\bbI_\ld \oplus \bbI_\ld \{1\} \oplus \bbI_\ld
    \{-1\}&\mbox{if }i\in I^-\mbox{ and }(\lambda_i, \lambda_{i+1})=(0,2),\\
    \mathbb{E}_{-i} \mathbb{E}_i\bbI_\ld \oplus \bbI_\ld&\mbox{if
    }(\af_i,\ld)=1,\\
    \bbE_{-i}\bbE_i\bbI_\ld&\mbox{if }(\af_i,\ld)=0;\end{cases}$$
\end{enumerate}
\end{prop}

Now we define the Huerfano-Khovanov 2-category $ \mathcal{HK}_{k,n} $
over the field $ \Bbbk $, $\mathrm{char}\Bbbk=2$.

\subsection{The objects}
The objects of $ \mathcal{HK}_{k,n} $ are the categories $
\mathcal{C}_{\lambda}$, $\ld\in P(V_{2\om_k})$.

\subsection{The 1-morphisms} For each $\ld\in P(2\om_k)$,
$\bbI_\ld\in\End_{\HK}(\ld)$ is the identity morphism and, for $\ld,\ld'\in P$,
set $\bbI_\ld\bbI_\ld'=\dt_{\ld,\ld'}\bbI_\ld$ as above. For each $i\in I$, we
have defined morphisms
$\bbE_i\bbI_\ld\in\Hom_{\HK}(\calC_\ld,\calC_{\ld+\af_i})$. Evidently, we have
$\bbE_i\bbI_\ld=\bbI_{\ld+\af_i}\bbE_i\bbI_\ld$. For $\ld,\ld'\in P(2\om_k)$,
we have
\[
\Hom_{\HK}(\calC_\ld,\calC_{\ld'})=\bigoplus_{\substack{\ui\in\Seq\\s\in\Z}}\C
\bbI_{\ld'}\bbE_{\ui}\bbI_{\ld}\{s\}
\]
where $\bbE_\ui:=\bbE_{i_1}\cdots \bbE_{i_r}\bbI_{\ld}$ if $\ui=(i_1,\ldots,i_r)\in I_\infty$,
and $s$ refers to a \emph{grading shift}. Observe that
$\bbI_{\ld'}\bbE_{\ui}\bbI_\ld=0$ unless $\cont(\ui)=\ld'-\ld$, and
$\bbI_{\ld+\cont(\ui)}\bbE_\ui\bbI_\ld=\bbE_\ui\bbI_\ld$.

\subsection{The 2-morphisms}
\label{2morphisms} Recall the convention
$(\ld_{-i},\ld_{-i+1})=(\ld_{i+1},\ld_i)$ for $i\in I^+$.

\begin{enumerate}
\item\textbf{The maps $ {1}_{i,\ld}$, ${1}_{\ld} $.}

Let $i\in I$, and let $ {1}_{i,\ld} \colon \mathbb{E}_i\bbI_\ld \rightarrow \mathbb{E}_i\bbI_\ld$,
and $ 1_{\ld} \colon \bbI_\ld \rightarrow \bbI_\ld $ be the identity maps.

\item\textbf{The maps $ y_{i;\ld}$.}

For $i\in I$ we define maps $ {y}_{i;\ld} \colon \mathbb{E}_i\bbI_\ld \rightarrow
\mathbb{E}_i\bbI_\ld $ of degree 2. Let $ T $ be the tangle diagram for the functor $
\mathbb{E}_i\bbI_\ld. $ It depends on the pair $ (\lambda_i, \lambda_{i+1}). $  Let $ a
$ and $ b $ be crossingless matches such that $ (Rb) T a $ is a disjoint union
of circles.  Thus $ \mathfrak{F}((Rb) T a) = (\A)^{\otimes p} $
for some natural number $ p. $ Define
$$ y_{i;\ld}((\beta_1 \otimes \cdots \otimes \beta_p)) = (\beta_1 \otimes \cdots \otimes x \beta_{j_i} \otimes \cdots \otimes \beta_p), $$ where
\begin{enumerate}
\item if $ (\lambda_i, \lambda_{i+1})=(1,2), $  then the $ j_i \text{th} $
factor in $(\A)^{\otimes p} $ corresponds to the circle passing
through the $ i\text{th} $ point on the bottom set of dots for tangle
$D_{\lambda, i} $ in Figure~\ref{d3}.
\item if $(\lambda_i, \lambda_{i+1})=(0,1), $  then the $ j_i \text{th} $
factor in $ (\A)^{\otimes p} $ corresponds to the circle passing
through the $ i\text{th} $ point on the top set of dots for tangle $D^{\lambda,
i} $ in Figure~\ref{d3}.
\item if $(\lambda_i, \lambda_{i+1})=(0,2), $  then the $ j_i \text{th} $
factor in $ (\A)^{\otimes p} $ corresponds to the circle passing
through the $ i\text{th} $ point on the top set of dots for tangle
$T^{\lambda,i} $ in Figure~\ref{d4}.
\item if $(\lambda_i, \lambda_{i+1})=(1,1), $  then the $ j_i \text{th} $
factor in $ (\A)^{\otimes p} $ corresponds to the circle passing
through the $ i\text{th} $ point on the bottom set of dots for tangle
$T_{\lambda,i} $ in Figure~\ref{d4}.
\end{enumerate}

\item\textbf{The map $ \cup_{i;\ld}$.}

We define a map $ {\cup}_{i;\ld} \colon \bbI_\ld \rightarrow \mathbb{E}_{-i}\mathbb{E}_i\bbI_\ld. $ There are four non-trivial cases for $ (\lambda_i,\lambda_{i+1}) $ to consider.
\begin{enumerate}
\item $ (\lambda_i, \lambda_{i+1})=(1,2). $  The identity functor is induced
from the identity tangle $ \I_{\lambda}$.  The functor $ \bbE_{-i} \mathbb{E}_i
$ is isomorphic to tensoring with the bimodule $ \mathfrak{F}(D^{\lambda+\alpha_i,i} \circ D_{\lambda,i}) $
which is equal to $ \mathfrak{F}(\I_{\lambda}) $.
Thus in this case $ {\cup}_{i;\ld} $ is given by the identity map.

\item $ (\lambda_i, \lambda_{i+1})=(1,1). $  Then the functor $
\mathbb{E}_{-i} \mathbb{E}_i $ is isomorphic to tensoring with the bimodule $ \mathfrak{F}(T^{\lambda+\alpha_i,i}
\circ T_{\lambda,i}) $.  Then $\cup_{i;\ld}$ is $
\mathfrak{F}(S^{\lambda,i}). $

\item $ (\lambda_i, \lambda_{i+1})=(0,2). $  Then the functor $
\mathbb{E}_{-i} \mathbb{E}_i $ is isomorphic to tensoring with the bimodule $
\mathfrak{F}(T_{\lambda+\alpha_i,i} \circ T^{\lambda, i}) =
\mathfrak{F}(\I_{\lambda}) \otimes \mathcal{A}. $ Then the bimodule map is
given by $ 1_{\ld} \otimes \iota. $

\item $ (\lambda_i, \lambda_{i+1})=(0,1). $ The functor $ \mathbb{E}_{-i}
\mathbb{E}_i $ is isomorphic to tensoring with the bimodule $ \mathfrak{F}(D_{\lambda +\alpha_i,i}
\circ D^{\lambda, i}). $ As in case 1, this tangle is isotopic to the identity
so the map between the functors is the identity map.
\end{enumerate}

\item\textbf{The map $ \cap_{i;\ld}$.}

We define a map $ {\cap}_{i;\ld} \colon \mathbb{E}_{-i} \mathbb{E}_i\bbI_\ld
\rightarrow \bbI_\ld. $ There are four non-trivial cases for $ (\lambda_i,
\lambda_{i+1}) $ to consider.
\begin{enumerate}
\item $ (\lambda_i, \lambda_{i+1})=(1,2). $ The functor $ \mathbb{E}_{-i}
\mathbb{E}_i $ is isomorphic to tensoring with the bimodule $ \mathfrak{F}(D^{\lambda+\alpha_i,i} \circ
D_{\lambda,i}) $ which is equal to $ \mathfrak{F}(\I_{\ld})$.
Thus in this case ${\cap}_{i;\ld} $ is given by the identity
map.
\item $ (\lambda_i, \lambda_{i+1})=(1,1). $ Then the functor $ \mathbb{E}_{-i}
\mathbb{E}_i $ is isomorphic to tensoring with the bimodule $ \mathfrak{F}(T^{\lambda+\alpha_i,i} \circ
T_{\lambda,i}) $. Then the homomorphism is $
\mathfrak{F}(S_{\lambda,i}). $
\item $ (\lambda_i, \lambda_{i+1})=(0,2). $  Then the functor $ \mathbb{E}_{-i}
\mathbb{E}_i $ is isomorphic to tensoring with the bimodule $ \mathfrak{F}(T_{\lambda+\alpha_i,i}
\circ T^{\lambda, i}) = \mathfrak{F}(\I_{\lambda}) \otimes \mathcal{A}. $
Then the bimodule map is given by $ 1_{\ld} \otimes \text{Tr}. $
\item $ (\lambda_i, \lambda_{i+1})=(0,1). $
The functor $ \mathbb{E}_{-i} \mathbb{E}_i $ is given by tensoring with the bimodule
$ \mathfrak{F}(D_{\lambda +\alpha_i,i} \circ D^{\lambda, i}). $
As in case 1, this tangle is isotopic to the identity so the map between the functors is the identity map.
\end{enumerate}

\item\textbf{The maps $ \psi_{i,j;\ld}$. }

We define a map $ {\psi}_{i,j;\ld} \colon \mathbb{E}_i\mathbb{E}_j\bbI_\ld
\rightarrow \mathbb{E}_j\mathbb{E}_i\bbI_\ld $ for $i,j\in I^{\pm}$.

There are four cases for $ i $ and $ j $ to consider and then subcases for $ \lambda. $
\begin{enumerate}
\item $ i=j. $  In this case, the functors are non-trivial only if $
\lambda_i = 0 $ and $ \lambda_{i+1} =2. $ The bimodule for $ \mathbb{E}_i
\mathbb{E}_i $ is isomorphic to tensoring with the bimodule $ \mathfrak{F}(T_{\lambda+\alpha_i,i} \circ
T^{\lambda,i}) = \mathfrak{F}(\I_{\lambda}) \otimes \mathcal{A}. $ Then $
\psi_{i,i} = 1_\ld \otimes \kappa. $

\item $ |i-j|>1. $  In this case, the functors $ \mathbb{E}_i \mathbb{E}_j $
and $ \mathbb{E}_j \mathbb{E}_i $ are isomorphic via an isomorphism induced
from a cobordism isotopic to the identity so set $ {\psi}_{i,j} $ to the
identity map.

\item $ {\psi}_{i,i+1} \colon \mathbb{E}_i \mathbb{E}_{i+1} \rightarrow
\mathbb{E}_{i+1} \mathbb{E}_i. $ There are four non-trivial subcases to
consider.
\begin{enumerate}
\item $ (\lambda_i, \lambda_{i+1}, \lambda_{i+2}) = (1,1,2). $
The bimodule for $ \mathbb{E}_i \mathbb{E}_{i+1} $ is $ \mathfrak{F}(D_{\lambda+\alpha_{i+1},i} \circ D_{\lambda,i+1}). $
The bimodule for $ \mathbb{E}_{i+1} \mathbb{E}_{i} $ is $ \mathfrak{F}(T_{\lambda+\alpha_{i},i+1} \circ T_{\lambda,i}). $
In this case we define the bimodule map to be
$ \mathfrak{F}(S^{\lambda,i,i+1}). $
\item $ (\lambda_i, \lambda_{i+1}, \lambda_{i+2}) = (1,1,1). $  The functor $ \mathbb{E}_i \mathbb{E}_{i+1} $
is given by tensoring with a bimodule isomorphic to
$$ \mathfrak{F}(D_{\lambda+\alpha_{i+1},i} \circ T_{\lambda,i+1}) \cong \mathfrak{F}(D_{\lambda+\alpha_{i+1},i} \circ T_{\lambda,i+1}) \otimes_{H_{\lambda}} \mathfrak{F}(\I_{\lambda}). $$
The bimodule for $ \mathbb{E}_{i+1} \mathbb{E}_i $ is isomorphic to $
\mathfrak{F}(D^{\lambda+\alpha_{i},i+1} \circ T_{\lambda,i}).  $ Then define $
{\psi}_{i,j} $ to be $ 1_\ld \otimes_{H_{\lambda}} \mathfrak{F}(S^{\lambda,i}) $ since
$$ \mathfrak{F}(D_{\lambda+\alpha_{i+1},i} \circ T_{\lambda,i+1}) \otimes_{H_{\lambda}} \mathfrak{F}(T^{\lambda+\alpha_i,-i} \circ T_{\lambda,i}) \cong
\mathfrak{F}(D^{\lambda+\alpha_{i},i+1} \circ T_{\lambda,i}).  $$

\item $ (\lambda_i, \lambda_{i+1}, \lambda_{i+2}) = (0,1,2). $
The bimodule for $ \mathbb{E}_i \mathbb{E}_{i+1} $ is isomorphic to
$$ \mathfrak{F}(T^{\lambda+\alpha_{i+1},i} \circ D_{\lambda,i+1}) \cong
\mathfrak{F}(\bbI_{\lambda+\alpha_i+\alpha_{i+1}})
\otimes_{H_{\lambda+\alpha_i+\alpha_{i+1}}}
\mathfrak{F}(T^{\lambda+\alpha_{i+1},i} \circ D_{\lambda,i+1}). $$ The bimodule
for $ \mathbb{E}_{i+1} \mathbb{E}_i $ is isomorphic to $
\mathfrak{F}(T^{\lambda+\alpha_{i},i+1} \circ D^{\lambda,i}). $ Then define $
{\psi}_{i,j} $ to be $\mathfrak{F}( S^{\lambda+\alpha_i+\alpha_{i+1},i} ) \otimes_{H_{\lambda}}
1_\ld $ since
$$ \mathfrak{F}(T^{\lambda+2\alpha_i+\alpha_{i+1},-(i+1)} \circ T_{\lambda+\alpha_i+\alpha_{i+1},i+1}) \otimes_{H_{\lambda+\alpha_i+\alpha_{i+1}}}
\mathfrak{F}(T^{\lambda+\alpha_{i+1},i} \circ D_{\lambda,i+1}) \cong \mathfrak{F}(T^{\lambda+\alpha_{i},i+1} \circ D^{\lambda,i}). $$

\item $ (\lambda_i, \lambda_{i+1}, \lambda_{i+2}) = (0,1,1). $  The bimodule
for $ \mathbb{E}_{i} \mathbb{E}_{i+1} $ is $
\mathfrak{F}(T^{\lambda+\alpha_{i+1},i} \circ T_{\lambda,i+1}). $ The bimodule
for $ \mathbb{E}_{i+1} \mathbb{E}_{i} $ is $
\mathfrak{F}(D^{\lambda+\alpha_i,i+1} \circ D^{\lambda,i}). $ Then set $
{\psi}_{i,j} = \mathfrak{F}(S_{\lambda,i+1,i}). $
\end{enumerate}

\item $ {\psi}_{i+1,i} \colon \mathbb{E}_{i+1} \mathbb{E}_{i} \rightarrow
\mathbb{E}_{i} \mathbb{E}_{i+1}. $

We essentially just have to read the maps in cases (c)(i)-(iv) above backwards.
\begin{enumerate}
\item $ (\lambda_i, \lambda_{i+1}, \lambda_{i+2}) = (1,1,2). $
The functors are just as in case (c)(i).  Now the map is
$ \mathfrak{F}(S_{\lambda,i,i+1}). $

\item $ (\lambda_i, \lambda_{i+1}, \lambda_{i+2}) = (1,1,1). $
The bimodule for $ \mathbb{E}_{i+1} \mathbb{E}_i $ is isomorphic to
$$ \mathfrak{F}(D^{\lambda+\alpha_{i},i+1} \circ T_{\lambda,i})  \cong \mathfrak{F}(D^{\lambda+\alpha_{i},i+1} \circ T_{\lambda,i}) \otimes_{H_{\lambda}} \mathfrak{F}(\I_{\lambda}). $$
Then define $ {\psi}_{i+1,i} = 1_\ld \otimes_{H_{\lambda}}
\mathfrak{F}(S^{\lambda,i+1}). $

\item $ (\lambda_i, \lambda_{i+1}, \lambda_{i+2}) = (0,1,2). $
The bimodule for $ \mathbb{E}_{i+1} \mathbb{E}_i $ is isomorphic to
$$ \mathfrak{F}(T^{\lambda+\alpha_{i},i+1} \circ D^{\lambda,i}) \cong \mathfrak{F}(\I_{\lambda+\alpha_i+\alpha_{i+1}}) \otimes_{H_{\lambda+\alpha_i + \alpha_{i+1}}}  \mathfrak{F}(T^{\lambda+\alpha_{i},i+1} \circ D^{\lambda,i}). $$
Then define $ {\psi}_{i+1,i} =
\mathfrak{F}(S^{\lambda+\alpha_i+\alpha_{i+1},i}) \otimes_{H_{\lambda}}
1_\ld . $

\item $ (\lambda_i, \lambda_{i+1}, \lambda_{i+2}) = (0,1,1). $
The functors are just as in case (c)(iv).  Now the map is
$ \mathfrak{F}(S^{\lambda,i+1,i}). $
\end{enumerate}
\end{enumerate}
\end{enumerate}

\begin{prop}
For all $i,j\in I$, and $ \ld \in P(V_{2\omega_k})$, the maps $ {y}_{i;\ld}, \psi_{i,j,\ld}, \cup_{i, \ld}, \cap_{i,\ld}$ are bimodule homomorphisms.
\end{prop}

For convenience of notation, we define the following 2-morphisms. If
$\theta\in\End(\bbE_\ui)$ let $ {\theta}^{[j]} = \underbrace{\theta \circ \cdots
\circ\theta}_j. $ For each $i\in I$, define the \emph{bubble}
$$\bubble^{N}_{i;\ld}=\cap_{i;\ld}\circ(1_{-i;\ld+\af_i}y_{i;\ld})^{[N]}\circ\cup_{i;\ld},$$
and define \emph{fake bubbles} inductively by the formula
\begin{align}\label{E:fake bubbles}
\left(\sum_{n\geq
0}\bubble_{i;\ld}^{(\af_{-i},\ld)-1+n}t^n\right)\left(\sum_{n\geq
0}\bubble_{-i;\ld}^{(\af_{i},\ld)-1+n}t^n\right)=1
\end{align}
and, $\bubble_{i;\ld}^{-1}=1$ whenever $(\af_i,\ld)=0$. Also, define \emph{half bubbles}
$$\ycup_{i;\ld}^N=(1_{-i;\ld+\af_i}y_{i;\ld})^{[N]} \circ
\cup_{i;\ld}\;\;\;\mbox{and}\;\;\;\ycap_{i;\ld}^N=\cap_{i;\ld} \circ (y_{i;\ld+\af_i} 1_{i,\ld})^{[N]}.$$
Finally, for $i,j\in I^\pm$, define
$$ \psi_{i,-j} = (1_{-j}1_i\cap_{-j}) \circ
(1_{-j}\psi_{j,i}1_{-j})\circ (\cup_j 1_i 1_{-j}). $$

\subsection{The 2-morphism relations}
Again, we will often omit the argument $ \ld $ when it is clear from context.

\subsection*{$\mathfrak{sl}_2$ relations.}

\begin{prop}
\label{kl3.1kl3.2}
For all $ i \in I$, $ ({\cap}_{-i}  {1}_i) \circ ({1}_i  {\cup}_i) = {1}_i = ({1}_i  {\cap}_i) \circ ({\cup}_{-i}  {1}_i). $
\end{prop}

\begin{proof}
The second equality is similar to the first equality.  The case $ i \in I^- $ is similar to the case $ i \in I^+$ so we just compute the map $ ({\cap}_i
 {1}_i) \circ ({1}_i  {\cup}_i)$ on the bimodule for the
functor $ \mathbb{E}_i$ for $ i \in I^+$. There are four cases to consider.

Suppose $ (\lambda_i, \lambda_{i+1})=(1,2). $  Then the tangle diagrams for the functors $ \mathbb{E}_i $ and $ \mathbb{E}_i \mathbb{E}_{-i} \mathbb{E}_i $ are
$ D_{\lambda,i} $ and $ D_{\lambda,i} \circ D^{\lambda+\alpha_i} \circ D_{\lambda,i} $
and can be found in Figure~\ref{d7}.

\begin{figure}[ht]
\[
\xy
(0,17)*{\mbox{\tiny2}};(10,17)*{\mbox{\tiny1}};
(0,10)*{};(10,20)*{}**\dir{-};(0,20)*{\bullet};(10,20)*{\bullet};(0,10)*{\bullet};
(10,-10)*{};
(10,10)*{\bullet};
(0,7)*{\mbox{\tiny1}};(10,7)*{\mbox{\tiny2}};
\endxy
\hspace{1in} \xy (0,17)*{\mbox{\tiny2}};(10,17)*{\mbox{\tiny1}};
(0,10)*{};(10,20)*{}**\dir{-};(0,20)*{\bullet};(10,20)*{\bullet};(0,20)*{\bullet};
(0,7)*{\mbox{\tiny1}};(10,7)*{\mbox{\tiny2}};
(0,10)*{};(10,0)*{}**\dir{-};(0,10)*{\bullet};(10,10)*{\bullet};(0,0)*{\bullet};
(0,-3)*{\mbox{\tiny2}};(10,-3)*{\mbox{\tiny1}};
(0,-10)*{};(10,0)*{}**\dir{-};(10,0)*{\bullet};(0,-10)*{\bullet};(10,-10)*{\bullet};
(0,-13)*{\mbox{\tiny1}};(10,-13)*{\mbox{\tiny2}};
\endxy
\]
\caption{Tangles for $\bbE_i$ and $\bbE_i\bbE_{-i}\bbE_i$,
$(\ld_i,\ld_{i+1})=(1,2)$, $i\in I^+$}\label{d7}
\end{figure}
The cobordism between the tangles is isotopic to the identity map so in this case the composition is equal to the identity map.

The case $ (\lambda_i, \lambda_{i+1})=(0,1) $ is similar to the $ (1,2) $ case.

Now let $ (\lambda_i, \lambda_{i+1})=(0,2). $
Then the tangle diagrams for the functors $ \mathbb{E}_i $ and $ \mathbb{E}_i \mathbb{E}_{-i} \mathbb{E}_i $ can be found in Figure~\ref{d8}.

\begin{figure}[ht]
\[
\xy (0,10); (0,7)*{\mbox{\tiny1}};(10,7)*{\mbox{\tiny1}};
(0,10)*{\bullet};(10,10)*{\bullet};
(0,10)*{};(10,10)*{}**\crv{(2,5)&(8,5)};
(0,0)*{\bullet};(10,0)*{\bullet};
(0,-3)*{\mbox{\tiny0}};(10,-3)*{\mbox{\tiny2}};
\endxy
\hspace{1in} \xy (0,7)*{\mbox{\tiny1}};(10,7)*{\mbox{\tiny1}};
(0,10)*{\bullet};(10,10)*{\bullet}; (0,10)*{};(10,10)*{}**\crv{(2,5)&(8,5)};
(0,-3)*{\mbox{\tiny0}};(10,-3)*{\mbox{\tiny2}};
(0,0)*{\bullet};(10,0)*{\bullet}; (5,-10)*\xycircle(5,5){-};
(0,-10)*{\bullet};(10,-10)*{\bullet};
(0,-13)*{\mbox{\tiny1}};(10.5,-13)*{\mbox{\tiny1}};
(0,-20)*{\bullet};(10,-20)*{\bullet};
(0,-23)*{\mbox{\tiny0}};(10,-23)*{\mbox{\tiny2}};
\endxy
\]
\caption{Tangles for $\bbE_i$ and $\bbE_i\bbE_{-i}\bbE_i$,
$(\ld_i,\ld_{i+1})=(0,2)$}\label{d8}
\end{figure}

Let $ B $ be the bimodule for the functor $ \mathbb{E}_i. $  Then the bimodule for $ \mathbb{E}_i \mathbb{E}_{-i} \mathbb{E}_i$ is isomorphic to $ \mathcal{A} \otimes B. $
The map $ \mathbb{E}_i \rightarrow \mathbb{E}_i \mathbb{E}_{-i} \mathbb{E}_i $ is given by the unit map which sends an element $ b \in B $ to $ 1 \otimes b.$
$ 1 \mapsto 1 \otimes b. $
The map $ \mathbb{E}_i \mathbb{E}_{-i} \mathbb{E}_i \rightarrow \mathbb{E}_i $ is obtained from the cobordism joining the circle to the upper cup which induces the multiplication map.
This maps $ 1 \otimes b $ to $ b. $
Thus the composition is equal to the identity.

Finally consider the case $ (\lambda_i, \lambda_{i+1})=(1,1). $
The tangle diagrams for the functors $ \mathbb{E}_i $ and $ \mathbb{E}_i \mathbb{E}_{-i} \mathbb{E}_i $ can be found in Figure~\ref{d9}.
\begin{figure}[ht]
\[
\xy (0,10); (0,7)*{\mbox{\tiny2}};(10,7)*{\mbox{\tiny0}};
(0,10)*{\bullet};(10,10)*{\bullet};
(0,0)*{};(10,0)*{}**\crv{(2,5)&(8,5)};
(0,0)*{\bullet};(10,0)*{\bullet};
(0,-3)*{\mbox{\tiny1}};(10,-3)*{\mbox{\tiny1}};
\endxy
\hspace{1in} \xy (0,7)*{\mbox{\tiny2}};(10,7)*{\mbox{\tiny0}};
(0,10)*{\bullet};(10,10)*{\bullet};
(0,-3)*{\mbox{\tiny1}};(10,-3)*{\mbox{\tiny1}};
(0,0)*{\bullet};(10,0)*{\bullet}; (5,0)*\xycircle(5,5){-};
(0,-10)*{\bullet};(10,-10)*{\bullet};
(0,-13)*{\mbox{\tiny2}};(10,-13)*{\mbox{\tiny0}};
(0,-10)*{\bullet};(10,-10)*{\bullet};
(0,-20)*{};(10,-20)*{}**\crv{(2,-15)&(8,-15)};
(0,-20)*{\bullet};(10,-20)*{\bullet};
(0,-23)*{\mbox{\tiny1}};(10,-23)*{\mbox{\tiny1}};
\endxy
\]
\caption{Tangles for $\bbE_i$ and $\bbE_i\bbE_{-i}\bbE_i$,
$(\ld_i,\ld_{i+1})=(1,1)$}\label{d9}
\end{figure}

Let $ B $ be the bimodule giving rise to the functor $ \mathbb{E}_i $ and $ \mathcal{A} \otimes B $ be the bimodule giving rise to the functor $ \mathbb{E}_i \bbE_{-i} \mathbb{E}_i. $
Let $ \alpha \otimes \beta \in B $ where $ \alpha $ is in the tensor factor corresponding to the circle passing through point $ i $ on the bottom row of the left side of Figure~\ref{d9} and $ \beta $ belongs to the remaining tensor factors.

The cobordism between the two tangle diagrams is a saddle which on the level of
bimodule maps, sends $ \alpha \otimes \beta \mapsto \Delta(\alpha) \otimes
\beta. $ Then the map from $ \mathbb{E}_i \bbE_{-i} \mathbb{E}_i$ to $
\mathbb{E}_i $ is given by $ \text{Tr} \otimes 1_\ld $ so $ \Delta(\alpha)
\otimes \beta \mapsto \alpha \otimes \beta $ by considering the two cases $
\alpha = 1 $ or $ x. $ Thus the composition is equal to the identity map.
\end{proof}

\begin{prop}
\label{kl3.3}
\begin{align*}
{y}_i = ({\cap}_{-i}  {1}_i) \circ ({1}_{i}  {y}_{-i}  {1}_{i}) \circ ({1}_i  {\cup}_i)= ({1}_i  {\cap}_i) \circ ({1}_i  {y}_{-i}  {1}_i) \circ ({\cup}_{-i} {1}_i).
\end{align*}
\end{prop}

\begin{proof}
We prove only the first equality as the second is similar.  There are four cases to consider for which the functor $ \bbE_{i} $ is non-zero.

Suppose $ (\lambda_i, \lambda_{i+1})=(1,2). $
Then the tangle diagrams for the functors $ \bbE_{i} $ and $ \bbE_{i} \mathbb{E}_{-i} \bbE_{i} $ can be found in Figure~\ref{d7}.

Let $ B $ be the bimodule for $ \bbE_{i} $ and $ \bbE_{i} \mathbb{E}_{-i}
\bbE_{i}. $  Let $ \alpha \otimes \beta \in B $ where $ \alpha $ is an
element in the tensor factor corresponding to a circle passing through point $
i $ in the bottom row of Figure~\ref{d7}. Then the first map $ {1}_i {\cup}_i
$ is given by the identity cobordism and is thus the identity.   The second map
is multiplication by $ x $ on all tensor components corresponding to circles
passing through the point $ i+1 $ in the second row of the right side of Figure
~\ref{d7}. The final map $ \bbE_{i} \mathbb{E}_{-i} \bbE_{i} \rightarrow
\bbE_{i} $ is also given by the identity cobordism.  Thus the composition
maps $ \alpha \otimes \beta \mapsto \alpha \otimes \beta \mapsto x\alpha
\otimes \beta \mapsto \alpha \otimes \beta. $ On the other hand, $ {y}_i
(\alpha \otimes \beta) = x \alpha \otimes \beta. $

The case $ (\lambda_i, \lambda_{i+1})=(0,1) $ is similar to the previous case.

Suppose $ (\lambda_i, \lambda_{i+1})=(0,2). $ Then the bimodule for the functor
$ \bbE_{i} $ is $ B = \mathfrak{F}(T^{\lambda,i}) $ and the tangle diagram
for $ \bbE_{i} \mathbb{E}_{-i} \bbE_{i} $  is $ \mathfrak{F}(T^{\lambda,i}
\circ T_{\lambda-\alpha_i,i} \circ T^{\lambda,i}) \cong \mathcal{A} \otimes B.
$ Let $ \alpha \otimes \beta \in B $ where $ \alpha $ is an element of the
tensor factor corresponding to the circle passing through the point $ i $ in
the top row of the tangle $ T^{\lambda,i} $ and $ \beta $ is an element in the
remaining tensor factors. Then the composition of maps send $ \alpha \otimes
\beta \mapsto 1 \otimes \alpha \otimes \beta \mapsto x \otimes \alpha \otimes
\beta \mapsto x\alpha \otimes \beta. $ This is equal to $ {y}_i(\alpha \otimes
\beta). $

Suppose $ (\lambda_i, \lambda_{i+1})=(1,1). $
Then the tangle diagrams for the functors $ \bbE_{i} $ and $ \bbE_{i} \mathbb{E}_{-i} \bbE_{i} $ can be found in Figure~\ref{d9}.

Let $ B $ be the bimodule for the functor $ \bbE_{-i} $ and $ \mathcal{A} \otimes B $ be the bimodule for $ \bbE_{i} \mathbb{E}_{-i} \bbE_{i}. $
Let $ \alpha \otimes \beta \in B $ where $ \alpha $ is an element in the tensor factor corresponding to the circle passing through point $ i $ on the bottom row of Figure ~\ref{d9}
and $ \beta $ is an element in the remaining tensor factors.
First let $ \alpha =1. $
Then
$$ 1 \otimes \beta \mapsto x \otimes 1 \otimes \beta + 1 \otimes x \otimes \beta \mapsto x \otimes x \otimes \beta \mapsto x \otimes \beta ={y}_i(1 \otimes \beta) $$
where the last map is $ \text{Tr} \otimes 1. $ If $ \alpha =x, $ then
$$ x \otimes \beta \mapsto x \otimes x \otimes \beta \mapsto 0 = {y}_i(x \otimes \beta). $$
\end{proof}

\begin{prop}
\label{kl3.4a}
Suppose $i\in I$ and $(-\af_i,\ld)>r+1$,  then $ \bubble_{i;\ld}^r = 0. $
\end{prop}

\begin{proof}
In order for $ r \geq 0, $ it must be the case that $ (-\af_i,\ld) \geq 2. $  Thus the only possibility is $ (\lambda_i, \lambda_{i+1}) = (0,2) $ and $ r =0. $
Then the bimodule for $ \mathbb{E}_{-i} \bbE_{i} $ is $ \mathcal{A} \otimes
\mathfrak{F}(\bbI_{\lambda}). $ Thus the map $ 1 \rightarrow \mathbb{E}_{-i}
\bbE_{i} $ is given by the unit map. The map $ \mathbb{E}_{-i} \bbE_{i}
\rightarrow 1 $ is given by the trace map. Thus the composition of the maps in
the proposition sends an element $ \beta\mapsto 1 \otimes \beta \mapsto
\text{Tr}(1) \otimes b = 0. $
\end{proof}

\begin{prop}
\label{kl3.4b}
If $ (\af_i,\ld)  \leq -1 $ then $ \bubble_{i;\ld}^{(-\af_i,\ld)-1} = 1. $
\end{prop}

\begin{proof}
The only cases to consider are $ (\lambda_i, \lambda_{i+1})=(0,2),(1,2),(0,1). $

Consider the case $ (0,2). $  Let $ B = \mathfrak{F}(\bbI_{\lambda}). $  Then the
bimodule corresponding to $ \mathbb{E}_{-i} \bbE_{i} $ is $ \mathcal{A}
\otimes B. $  Let $ \beta \in B. $ Then $ {\cup}_i (\beta) = 1 \otimes \beta$,
$ {y}_i(1 \otimes \beta)=x \otimes \beta$,  and $ {\cap}_i(x \otimes
\beta) = \text{Tr}(x)\beta = \beta. $ Thus in this case, the composition is the
identity map.

For the case $ (1,2), $
$ (-\af_i,\ld)-1=0. $
The cobordism between the tangle diagrams for the identity functor and $ \mathbb{E}_{-i} \bbE_{i} $ is isotopic to the identity cobordism.
Similarly, the cobordism between the tangle diagrams for the functors $ \mathbb{E}_{-i} \bbE_{i} $ and the identity functor is isotopic to the identity cobordism.  Thus the bimodule map is equal to the
identity.

The case $ (0,1) $ is the same as the case $ (1,2). $
\end{proof}

\begin{prop}
\label{kl3.6}
Let $i\in I$. If $(\af_i,\ld)\geq 1$, then
\begin{align*}
1_{i;\ld-\af_i}  1_{-i;\ld} =  {\psi}_{-i,i;\ld} \circ {\psi}_{i,-i;\ld}+
\sum_{f=0}^{(\af_i,\ld)-1}
\sum_{g=0}^f \ycup_{-i;\ld}^{(\af_i,\ld)-f-1} \circ \bubble^{-(\af_i,\ld)-1+g}_{i;\ld}
 \circ \ycap_{-i;\ld}^{f-g}
\end{align*}
\end{prop}

\begin{proof}
There are three cases to consider: $ (\lambda_i, \lambda_{i+1})=(1,0), (2,1), (2,0). $

For the case $ (1,0), $ the first term on the right hand side is zero since that map passes through the functor $ \mathbb{E}_i \mathbb{E}_i \bbE_{-i} $ which is zero for this $ \lambda. $
The summation on the right hand side reduces to
$$ {\ycup}_{-i;\ld}^{0} \circ \bubble^{-2}_{i;\ld} \circ {\ycap}_{i;\ld}^0 =
{\cup}_{-i;\ld}  \circ {\cap}_{-i; \ld} $$ by definition \eqref{E:fake bubbles} of the fake bubbles. This map is a
composition $ \mathbb{E}_i \bbE_{-i} \rightarrow 1 \rightarrow \mathbb{E}_i
\bbE_{-i}. $ This composition of maps is the identity.

The case $ (2,1) $ is similar to the $ (1,0) $ case.

For the case $ (2,0), $ the first term on the right hand side is zero as in the previous two cases.
The summation on the right hand side consists of three terms which simplifies by \eqref{E:fake bubbles} to:
$$ \ycup_{-i;\ld}^1 \circ \cap_{-i;\ld} + {\cup}_{-i;\ld} \circ \ycap_{-i;\ld}^1  +\cup_{-i;\ld}\circ\bubble_{i;\ld}^2\circ {\cap}_{-i;\ld}. $$
Let $ B = \mathfrak{F}(\bbI_{\lambda}). $  Then the bimodule for $ \mathbb{E}_i
\bbE_{-i} $ is $ \mathcal{A} \otimes B. $ Then
$$ \ycup_{-i;\ld}^1 \circ {\cap}_{-i;\ld} \colon \mathbb{E}_i \bbE_{-i} \rightarrow \bbI \rightarrow \mathbb{E}_i \bbE_{-i} \rightarrow \mathbb{E}_i \bbE_{-i}. $$
Under this composition of maps, $ 1 \otimes b $ maps to zero since the first map is given by a trace map on the first component.  The element $ x \otimes b $ gets mapped to $ x \otimes b $ as follows:
$$ x \otimes b \mapsto b \mapsto 1 \otimes b \mapsto x \otimes b, $$
where the first map is the trace map, the second map is the unit map and the third map is multiplication by $ x. $
Similarly,
$$ {\cup}_{-i;\ld} \circ \ycap_{-i;\ld}^1 \colon \mathbb{E}_i \bbE_{-i} \rightarrow \mathbb{E}_i \bbE_{-i}\rightarrow \bbI  \rightarrow \mathbb{E}_i \bbE_{-i}. $$
Under this composition, $ 1 \otimes b \mapsto 1 \otimes b $ and $ x \otimes b \mapsto 0. $
Finally, the map
$$ \cup_{-i;\ld}\circ\bubble_{i;\ld}^2\circ {\cap}_{-i;\ld}$$ is zero because the middle term is zero.
Thus the right hand side is the identity as well.
\end{proof}

\begin{prop}
\label{kl3.5}
Let $i\in I^+$.
\begin{enumerate}
\item If $(\af_{i},\ld)\leq0$, then
$$({1}_i  {\cap}_{-i;\ld}) \circ ({\psi}_{i,i;\ld-\alpha_i}  {1}_{-i}) \circ ({1}_i
 {\cup}_{-i;\ld}) =  \sum_{f=0}^{-(\af_i,\ld)}
{{y}_i}^{-(\af_i,\ld)-f} \bubble_{-i;\ld}^{(\af_i,\ld)-1+f}. $$

\item If $(\af_i,\ld)\geq-2$, then
$$({\cap}_{i;\ld+\alpha_i}  {1}_i) \circ ({1}_i  {\psi}_{i,i;\ld}) \circ ({\cup}_{i;\ld+\af_i}
 {1}_i) = \sum_{g=0}^{(\af_i,\ld)+2} \bubble_{i,\ld}^{-(\af_{i},\ld)-3+g}
{{y}_i}^{(\af_i,\ld)-g+2}. $$
\end{enumerate}
\end{prop}

\begin{proof} We prove (1), the proof of (2) being similar.
Since the map on both sides pass through the functor $ \mathbb{E}_i \mathbb{E}_i \bbE_{-i}, $ the maps on both sides are zero unless $ (\lambda_i, \lambda_{i+1})=(1,1). $
The functors for $ \mathbb{E}_i $ and $ \mathbb{E}_i \mathbb{E}_i \bbE_{-i} $ are given by tangles in Figure~\ref{d9}.

Let $ B $ be the bimodule for the functor $ \mathbb{E}_i $ so $ \mathcal{A} \otimes B $ is the bimodule for the functor $ \mathbb{E}_i \mathbb{E}_i \bbE_{-i}. $
Let $ \alpha \otimes \beta \in B $ where $ \alpha $ is an element in the tensor factor corresponding to a circle passing through point $ i $ in the bottom row of the left side of figure ~\ref{d9} and
$ \beta $ is an element in the other tensor factors.
Consider first $ \alpha =1. $
The left hand side maps an element $ \alpha \otimes \beta $ as follows:
$$ 1 \otimes \beta \mapsto x \otimes 1 \otimes \beta + 1 \otimes x \otimes \beta \mapsto 1 \otimes 1 \otimes \beta \mapsto 1 \otimes \beta $$
where the first map is $ \Delta \otimes 1, $ the second map is $ \kappa \otimes
1 \otimes 1 $ and the third map is $ m \otimes 1. $ If $ \alpha = x, $ the left
hand maps $ \alpha \otimes \beta $ as follows:
$$ x \otimes \beta \mapsto x \otimes x \otimes \beta \mapsto 1 \otimes x \otimes \beta \mapsto x \otimes \beta. $$
The right hand side is $ - 1 $ by convention.
\end{proof}

\subsection*{nil-Hecke relations.}

\begin{prop}
\label{kl3.8a}  For $ i \in I^+, $ $ {\psi}_{i,i}^{[2]} = 0. $
\end{prop}

\begin{proof}
Since $ \mathbb{E}_i \mathbb{E}_i $ is identically zero unless $ (\lambda_i,
\lambda_{i+1}) = (0,2), $ we need only consider this case. Let $ B =
\mathfrak{F}(\bbI_{\lambda}). $ Then the bimodule for $ \mathbb{E}_i \mathbb{E}_i
$ is isomorphic to $ \mathfrak{F}(T_{\lambda,i} \circ T^{\lambda,i}) = \mathcal{A}
\otimes B. $

Then $ {\psi}_{i,i} \circ {\psi}_{i,i} \colon \mathcal{A} \otimes B \rightarrow
\mathcal{A} \otimes B \rightarrow \mathcal{A} \otimes B. $ This map sends $ 1
\otimes b \mapsto 0 $ and $ x \otimes b \mapsto 1 \otimes b \mapsto 0. $
\end{proof}

\begin{prop}
\label{kl3.8b} Let $ i \in I^+. $  Then, $ ({\psi}_{i,i}  {1}_i) \circ ({1}_i  {\psi}_{i,i}) \circ
({\psi}_{i,i}  {1}_i) = ({1}_i  {\psi}_{i,i}) \circ ({\psi}_{i,i} {1}_i) \circ
({1}_i  {\psi}_{i,i}). $
\end{prop}

\begin{proof}
Both sides are natural transformations of the functor $ \mathbb{E}_i\mathbb{E}_i\mathbb{E}_i. $  However, by definition this composition is zero.
\end{proof}

\begin{prop}
\label{kl3.9}
For $ i \in I^+, $
$ ({1}_i  {1}_i) = ({\psi}_{i,i}) \circ ({y}_i  {1}_i) - ({1}_i
 {y}_i) \circ ({\psi}_{i,i}) = ({y}_i  {1}_i) \circ
({\psi}_{i,i}) - ({\psi}_{i,i}) \circ ({1}_i  {y}_i). $
\end{prop}

\begin{proof}
The only case to check is
$ (\lambda_i, \lambda_{i+1})=(0,2) $
since otherwise $ \mathbb{E}_i \mathbb{E}_i = 0. $
Let $ B = \mathfrak{F}(\bbI_{\lambda}). $ Then
the bimodule for $ \mathbb{E}_i \mathbb{E}_i $ is isomorphic to $ \mathcal{A} \otimes B. $
Then
$$ ({\psi}_{i,i}) \circ ({y}_i  {1}_i) \colon \mathcal{A} \otimes B\rightarrow \mathcal{A} \otimes B. $$
Under this map, $ 1 \otimes b \mapsto x \otimes b \mapsto 1 \otimes b $ and $ x
\otimes b \mapsto 0. $ For the map $ ({1}_i  {y}_i) \circ ({\psi}_{i,i}), $ $ 1
\otimes b \mapsto 0, $ and $ x \otimes b \mapsto 1 \otimes b \mapsto x \otimes
b. $ This gives the first equality since our field is $ \mathbb{F}_2.$

For the second equality, $ ({y}_i  {1}_i) \circ ({\psi}_{i,i})\colon 1 \otimes
b \mapsto 0, $ $ ({y}_i  {1}_i) \circ ({\psi}_{i,i}) \colon x \otimes b \mapsto
1 \otimes b \mapsto x \otimes b. $ Similarly, $ ({\psi}_{i,i}) \circ ({1}_i
{y}_i) \colon 1 \otimes b \mapsto x \otimes b \mapsto 1 \otimes b $ and $
({\psi}_{i,i}) \circ ({1}_i  {y}_i) \colon x \otimes b \mapsto 0. $
\end{proof}

\begin{prop}
\label{kl3.10} For $i,j\in I^-$,

\begin{align*}
{\psi}_{j,i} = &({\cap}_{-j}  {1}_i  {1}_j) \circ ({1}_j  {\cap}_{-i}  {1}_{-j}
 {1}_i  {1}_j) \circ
({1}_j  {1}_i  {\psi}_{-j,-i}  {1}_i  {1}_j) \circ
({1}_j  {1}_i  {1}_{-j}  {\cup}_i  {1}_j) \circ ({1}_j  {1}_i  {\cup}_j)\\
=&({1}_i  {1}_j  {\cap}_i) \circ ({1}_i  {1}_j  {1}_{-i} {\cap}_j
 {1}_i) \circ
({1}_i  {1}_j  {\psi}_{-j,-i}  {1}_j  {1}_i) \circ
({1}_i  {\cup}_{-j}  {1}_{-i}  {1}_j  {1}_i) \circ ({\cup}_{-i}  {1}_j {1}_i).
\end{align*}
\end{prop}

\begin{proof} Let $i,j\in I^-$.
We prove only the first equality. If $ |i-j|>1, $ the proposition is easy because then $ {\psi}_{\pm i,\pm j}$ are identity morphisms. Therefore, we take $ i = j+1$, the case $ i =j-1 $ being similar. The natural transformation on the right side of the proposition is a composition of natural transformations:
$$ \bbE_{j} \bbE_{j+1} \rightarrow
 \bbE_{j} \bbE_{j+1} \mathbb{E}_{-j} \bbE_{j} \rightarrow
 \bbE_{j} \bbE_{j+1} \mathbb{E}_{-j} \mathbb{E}_{-j-1}  \bbE_{j+1} \bbE_{j} \rightarrow
 \bbE_{j} \bbE_{j+1} \mathbb{E}_{-j-1} \mathbb{E}_{-j}  \bbE_{j+1} \bbE_{j} \rightarrow
 \bbE_{j} \mathbb{E}_{-j} \bbE_{j+1} \bbE_{j} \rightarrow
 \bbE_{j+1} \bbE_{j}. $$

There are four nontrivial cases for $ \lambda. $ We prove the case $ (\lambda_j, \lambda_{j+1}, \lambda_{j+2})=(2,1,1)$. The proofs of the remaining cases $(2,1,0)$, $(1,1,0)$, and $(1,1,1)$ are similar.

\begin{figure}[ht]
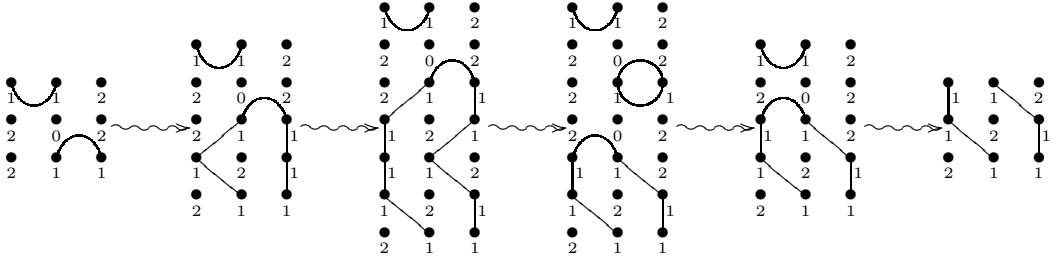

\[
\xy (0,0)*{
\xy
(0,10)*{\bullet};(6,10)*{\bullet};(12,10)*{\bullet};(0,10)*{};(6,10)*{}**\crv{(1,6)&(5,6)};
(0,8)*{\mbox{\tiny1}};(6,8)*{\mbox{\tiny1}};(12,8)*{\mbox{\tiny2}};
(0,5)*{\bullet};(6,5)*{\bullet};(12,5)*{\bullet};
(0,3)*{\mbox{\tiny2}};(6,3)*{\mbox{\tiny0}};(12,3)*{\mbox{\tiny2}};
(0,0)*{\bullet};(6,0)*{\bullet};(12,0)*{\bullet};(6,0)*{};(12,0)*{}**\crv{(7,4)&(11,4)};
(0,-2)*{\mbox{\tiny2}};(6,-2)*{\mbox{\tiny1}};(12,-2)*{\mbox{\tiny1}};
\endxy}="a";
(25,0)*{
\xy
(0,20)*{\bullet};(6,20)*{\bullet};(12,20)*{\bullet};(0,20)*{};(6,20)*{}**\crv{(1,16)&(5,16)};
(0,18)*{\mbox{\tiny1}};(6,18)*{\mbox{\tiny1}};(12,18)*{\mbox{\tiny2}};
(0,15)*{\bullet};(6,15)*{\bullet};(12,15)*{\bullet};
(0,13)*{\mbox{\tiny2}};(6,13)*{\mbox{\tiny0}};(12,13)*{\mbox{\tiny2}};
(0,10)*{\bullet};(6,10)*{\bullet};(12,10)*{\bullet};(6,10)*{};(12,10)*{}**\crv{(7,14)&(11,14)};
(0,8)*{\mbox{\tiny2}};(6,8)*{\mbox{\tiny1}};(13,8)*{\mbox{\tiny1}};
(0,5)*{\bullet};(6,5)*{\bullet};(12,5)*{\bullet};(0,5)*{};(6,10)*{}**\dir{-};
(0,3)*{\mbox{\tiny1}};(6,3)*{\mbox{\tiny2}};(13,3)*{\mbox{\tiny1}};
(0,0)*{\bullet};(6,0)*{\bullet};(12,0)*{\bullet};(6,0)*{};(0,5)*{}**\dir{-};(12,0)*{};(12,10)*{}**\dir{-};
(0,-2)*{\mbox{\tiny2}};(6,-2)*{\mbox{\tiny1}};(12,-2)*{\mbox{\tiny1}};
\endxy}="b";
(50,0)*{
\xy
(0,30)*{\bullet};(6,30)*{\bullet};(12,30)*{\bullet};(0,30)*{};(6,30)*{}**\crv{(1,26)&(5,26)};
(0,28)*{\mbox{\tiny1}};(6,28)*{\mbox{\tiny1}};(12,28)*{\mbox{\tiny2}};
(0,25)*{\bullet};(6,25)*{\bullet};(12,25)*{\bullet};
(0,23)*{\mbox{\tiny2}};(6,23)*{\mbox{\tiny0}};(12,23)*{\mbox{\tiny2}};
(0,20)*{\bullet};(6,20)*{\bullet};(12,20)*{\bullet};(6,20)*{};(12,20)*{}**\crv{(7,24)&(11,24)};
(0,18)*{\mbox{\tiny2}};(6,18)*{\mbox{\tiny1}};(13,18)*{\mbox{\tiny1}};
(0,15)*{\bullet};(6,15)*{\bullet};(12,15)*{\bullet};(0,15)*{};(6,20)*{}**\dir{-};(12,15)*{};(12,20)*{}**\dir{-};
(1,13)*{\mbox{\tiny1}};(6,13)*{\mbox{\tiny2}};(12,13)*{\mbox{\tiny1}};
(0,10)*{\bullet};(6,10)*{\bullet};(12,10)*{\bullet};(6,10)*{};(12,15)*{}**\dir{-};
(1,8)*{\mbox{\tiny1}};(6,8)*{\mbox{\tiny1}};(12,8)*{\mbox{\tiny2}};
(0,5)*{\bullet};(6,5)*{\bullet};(12,5)*{\bullet};(0,5)*{};(0,15)*{}**\dir{-};(12,5)*{};(6,10)*{}**\dir{-};
(0,3)*{\mbox{\tiny1}};(6,3)*{\mbox{\tiny2}};(13,3)*{\mbox{\tiny1}};
(0,0)*{\bullet};(6,0)*{\bullet};(12,0)*{\bullet};(6,0)*{};(0,5)*{}**\dir{-};(12,0)*{};(12,5)*{}**\dir{-};
(0,-2)*{\mbox{\tiny2}};(6,-2)*{\mbox{\tiny1}};(12,-2)*{\mbox{\tiny1}};
\endxy}="c";
(75,0)*{
\xy
(0,30)*{\bullet};(6,30)*{\bullet};(12,30)*{\bullet};(0,30)*{};(6,30)*{}**\crv{(1,26)&(5,26)};
(0,28)*{\mbox{\tiny1}};(6,28)*{\mbox{\tiny1}};(12,28)*{\mbox{\tiny2}};
(0,25)*{\bullet};(6,25)*{\bullet};(12,25)*{\bullet};
(0,23)*{\mbox{\tiny2}};(6,23)*{\mbox{\tiny0}};(12,23)*{\mbox{\tiny2}};
(0,20)*{\bullet};(6,20)*{\bullet};(12,20)*{\bullet};(9,20)*\xycircle(3,3){-};
(0,18)*{\mbox{\tiny2}};(6,18)*{\mbox{\tiny1}};(13,18)*{\mbox{\tiny1}};
(0,15)*{\bullet};(6,15)*{\bullet};(12,15)*{\bullet};
(0,13)*{\mbox{\tiny2}};(6,13)*{\mbox{\tiny0}};(12,13)*{\mbox{\tiny2}};
(0,10)*{\bullet};(6,10)*{\bullet};(12,10)*{\bullet};(0,10)*{};(6,10)*{}**\crv{(1,14)&(5,14)};
(1,8)*{\mbox{\tiny1}};(6,8)*{\mbox{\tiny1}};(12,8)*{\mbox{\tiny2}};
(0,5)*{\bullet};(6,5)*{\bullet};(12,5)*{\bullet};(0,5)*{};(0,10)*{}**\dir{-};(12,5)*{};(6,10)*{}**\dir{-};
(0,3)*{\mbox{\tiny1}};(6,3)*{\mbox{\tiny2}};(13,3)*{\mbox{\tiny1}};
(0,0)*{\bullet};(6,0)*{\bullet};(12,0)*{\bullet};(6,0)*{};(0,5)*{}**\dir{-};(12,0)*{};(12,5)*{}**\dir{-};
(0,-2)*{\mbox{\tiny2}};(6,-2)*{\mbox{\tiny1}};(12,-2)*{\mbox{\tiny1}};
\endxy}="d";
(100,0)*{
\xy
(0,20)*{\bullet};(6,20)*{\bullet};(12,20)*{\bullet};(0,20)*{};(6,20)*{}**\crv{(1,16)&(5,16)};
(0,18)*{\mbox{\tiny1}};(6,18)*{\mbox{\tiny1}};(12,18)*{\mbox{\tiny2}};
(0,15)*{\bullet};(6,15)*{\bullet};(12,15)*{\bullet};
(0,13)*{\mbox{\tiny2}};(6,13)*{\mbox{\tiny0}};(12,13)*{\mbox{\tiny2}};
(0,10)*{\bullet};(6,10)*{\bullet};(12,10)*{\bullet};(0,10)*{};(6,10)*{}**\crv{(1,14)&(5,14)};
(1,8)*{\mbox{\tiny1}};(6,8)*{\mbox{\tiny1}};(12,8)*{\mbox{\tiny2}};
(0,5)*{\bullet};(6,5)*{\bullet};(12,5)*{\bullet};(0,5)*{};(0,10)*{}**\dir{-};(12,5)*{};(6,10)*{}**\dir{-};
(0,3)*{\mbox{\tiny1}};(6,3)*{\mbox{\tiny2}};(13,3)*{\mbox{\tiny1}};
(0,0)*{\bullet};(6,0)*{\bullet};(12,0)*{\bullet};(6,0)*{};(0,5)*{}**\dir{-};(12,0)*{};(12,5)*{}**\dir{-};
(0,-2)*{\mbox{\tiny2}};(6,-2)*{\mbox{\tiny1}};(12,-2)*{\mbox{\tiny1}};
\endxy}="e";
(125,0)*{
\xy
(0,10)*{\bullet};(6,10)*{\bullet};(12,10)*{\bullet};
(1,8)*{\mbox{\tiny1}};(6,8)*{\mbox{\tiny1}};(12,8)*{\mbox{\tiny2}};
(0,5)*{\bullet};(6,5)*{\bullet};(12,5)*{\bullet};(0,5)*{};(0,10)*{}**\dir{-};(12,5)*{};(6,10)*{}**\dir{-};
(0,3)*{\mbox{\tiny1}};(6,3)*{\mbox{\tiny2}};(13,3)*{\mbox{\tiny1}};
(0,0)*{\bullet};(6,0)*{\bullet};(12,0)*{\bullet};(6,0)*{};(0,5)*{}**\dir{-};(12,0)*{};(12,5)*{}**\dir{-};
(0,-2)*{\mbox{\tiny2}};(6,-2)*{\mbox{\tiny1}};(12,-2)*{\mbox{\tiny1}};
\endxy}="f";
{\ar@{~>}"a";"b"};{\ar@{~>}"b";"c"};{\ar@{~>}"c";"d"};{\ar@{~>}"d";"e"};{\ar@{~>}"e";"f"};
\endxy
\]
\caption{Tangles for compositions of natural transformations in the $(2,1,1)$ case.}\label{d16}
\end{figure}

Let $ B $ be the bimodule representing the functor $ \bbE_j \bbE_{j+1} $ and $ B' $ the bimodule representing the functor $ \bbE_{j+1} \bbE_j. $
Then the morphism is the composition $ B \rightarrow B \rightarrow B \rightarrow \mathcal{A} \otimes B \rightarrow B \rightarrow B' $
induced by the tangle cobordisms in Figure ~\ref{d16}. The first and second
maps are the identity maps.  The third map is comultiplication.  The fourth map
is the counit map and the last map is $ {\psi}_{j,j+1}. $ Computing this
composition on elements as in previous propositions easily gives that it is
equal to $ {\psi}_{j,j+1}. $
\end{proof}

\subsection*{$R(\nu)$ relations.}

\begin{prop}
\label{kl3.13} For $i,j\in I^{\pm}$, $i\neq j$,
 $$ {\psi}_{-j,i} \circ {\psi}_{i,-j} = {1}_i  {1}_{-j}. $$
\end{prop}

\begin{proof} Note that for $ |i-j|>1, $ the left hand side is easily seen to be the identity so let $ j=i+1. $  The case $ j=i-1 $ is similar.  Thus the left hand side is:
\begin{align*}
{\psi}_{-j,i}\circ
{\psi}_{i,-j} \colon\mathbb{E}_i \bbE_{-i-1} \rightarrow
\mathbb{E}_{-i-1} \mathbb{E}_{i+1} \mathbb{E}_{i} \mathbb{E}_{-i-1} &\rightarrow
\mathbb{E}_{-i-1} \mathbb{E}_{i} \mathbb{E}_{i+1} \mathbb{E}_{-i-1} \rightarrow
\mathbb{E}_{-i-1} \mathbb{E}_{i}\rightarrow\\ &\rightarrow
\mathbb{E}_{-i-1} \mathbb{E}_{i} \mathbb{E}_{i+1} \mathbb{E}_{-i-1} \rightarrow
\mathbb{E}_{-i-1} \mathbb{E}_{i+1} \mathbb{E}_{i} \mathbb{E}_{-i-1} \rightarrow
\mathbb{E}_i \mathbb{E}_{-i-1}.
\end{align*}

There are four non-trivial cases for $ \lambda. $

Case 1: $ (\lambda_i, \lambda_{i+1}, \lambda_{i+2})=(1,2,1). $
Let $ B $ be the bimodule representing the functor
$ \mathbb{E}_i \mathbb{E}_{-i-1}. $
Then
$$ {\psi}_{-j,i} \circ{\psi}_{i,-j} \colon B \rightarrow \mathcal{A} \otimes B \rightarrow B \rightarrow B \rightarrow B \rightarrow \mathcal{A} \otimes B \rightarrow B. $$
The first map is $ \iota \otimes 1_\ld. $  The second map is multiplication
$ m. $  The third and fourth maps are the identity.  The fifth map is
comultiplication $ \Delta. $  The last map is $ \text{Tr} \otimes 1. $ It
is easy to check on elements that this is the identity map.

Case 2: $ (\lambda_i, \lambda_{i+1}, \lambda_{i+2})=(1,2,0). $
Let $ B $ be the bimodule representing the functor
$ \mathbb{E}_i \mathbb{E}_{-i-1}. $
Then
$$ {\psi}_{-j,i} \circ {\psi}_{i,-j} \colon
B \rightarrow B \rightarrow \mathcal{A} \otimes B  \rightarrow B \rightarrow
\mathcal{A} \otimes B \rightarrow B \rightarrow B. $$ The first map is the
identity.  The second map is $ \Delta $ by lemma ~\ref{mulemma1}.  The third map is $ \text{Tr} \otimes
1$ where the trace map is applied to the tensor factor arising from the
new circle component. The fourth map is $ \iota \otimes 1. $  The fifth
map is multiplication by lemma ~\ref{mulemma2}.  The last map is the identity. It is easy to check that
this composition is the identity on all elements.

Case 3: $ (\lambda_i, \lambda_{i+1}, \lambda_{i+2})=(0,2,1). $
This is similar to case 2.

Case 4: $ (\lambda_i, \lambda_{i+1}, \lambda_{i+2})=(0,2,0). $
This is similar to case 1.
\end{proof}

\begin{prop}
\label{kl3.14a} If $i,j\in I^+$ and $|i-j|>1, $ then $ {\psi}_{j,i} \circ {\psi}_{i,j} = {1}_i
{1}_j. $
\end{prop}

\begin{proof}
The tangle diagrams for the bimodules for $ \mathbb{E}_i \mathbb{E}_j $ and $ \mathbb{E}_j \mathbb{E}_i $ are the same up to isotopy.  The maps in the proposition are obtained from cobordisms isotopic to the identity so they are identity maps.
\end{proof}

\begin{prop}
\label{kl3.14b} If $i,j\in I^+$ and $|i-j|=1$ then $ {\psi}_{j,i} \circ {\psi}_{i,j} =
({y}_i  {1}_j + {1}_i  {y}_j). $
\end{prop}

\begin{proof}
Assume $ j=i+1. $  The case $ j = i-1 $ is similar.  There are eight cases for $ \lambda $ such that $ \mathbb{E}_i \mathbb{E}_{i+1} $ is non-zero.
In all cases let $ a $ and $ b $ be cup diagrams.  Let $ B $ be the bimodule for $ \mathbb{E}_i \mathbb{E}_{i+1} $ and $ B' $ the bimodule for $ \mathbb{E}_{i+1} \mathbb{E}_i. $

Case 1: $ (\lambda_i, \lambda_{i+1}, \lambda_{i+2})=(0,0,1). $ Since $
\mathbb{E}_{i+1} \mathbb{E}_{i}=0, $  the map $ {\psi}_{i+1,i} \circ
{\psi}_{i,i+1} = 0. $  The bimodule representing the functor $ \mathbb{E}_i
\mathbb{E}_{i+1} $ is isomorphic to $ \mathfrak{F}(D^{\lambda+\alpha_{i+1},i} \circ
D^{\lambda,i+1}). $ Since the circle passing through point $ i $ on the bottom
row of $ D^{\lambda+\alpha_{i+1},i} \circ D^{\lambda,i+1} $ is the same as the
circle passing through point $ i+1 $ in the middle row, the map on the right
side of the proposition is zero as well.

Case 2: $ (\lambda_i, \lambda_{i+1}, \lambda_{i+2})=(1,0,1). $
This is similar to case 1.

Case 3: $ (\lambda_i, \lambda_{i+1}, \lambda_{i+2})=(1,0,2). $
This is similar to case 1.

Case 4: $ (\lambda_i, \lambda_{i+1}, \lambda_{i+2})=(0,0,2). $
This is similar to case 1.

Case 5: $ (\lambda_i, \lambda_{i+1}, \lambda_{i+2})=(0,1,1). $  In this case $ B
\cong \mathfrak{F}(T^{\lambda+\alpha_{i+1},i} \circ T_{\lambda, i+1}) $ and $ B' \cong
\mathfrak{F}(D^{\lambda+\alpha_i,i+1} \circ D^{\lambda,i}). $
Let $ a $ and $ b $ be crossingless matches.
\begin{itemize}
\item
Suppose that the
circle passing through point $ i+1 $ on the bottom row of $
{}_a(T^{\lambda+\alpha_{i+1},i}) \circ T_{\lambda, i+1})_b $ is the same as the circle
passing through point $ i $ of the top row. Then $ {}_a B_{b} = \mathcal{A}
\otimes R $ and $ {}_a B_{b}' = \mathcal{A} \otimes \mathcal{A} \otimes R $
where $ R $ is a tensor product of $ \mathcal{A} $ corresponding to the
remaining circles. Then the map on the left side of the proposition is $ (m
\otimes 1) \circ (\Delta \otimes 1). $ Thus it maps an element $ 1 \otimes r $
to $ 2x \otimes r. $ On the other hand, $ {y}_i(1 \otimes r) = +x \otimes r. $
Also, $ {y}_{i+1}(1 \otimes r) = x \otimes r. $  Thus both sides are the same.

\item Suppose that the circle passing through point $ i+1 $ on the bottom is
different from the circle passing through point $ i $ on the top. Then $ {}_a
B_{b} = \mathcal{A} \otimes \mathcal{A} \otimes R $ and $ {}_a B_{b}' =
\mathcal{A} \otimes R. $ Then the map on the left side of the proposition is $
(\Delta \otimes 1_\ld) \circ (m \otimes 1_\ld). $ Thus it maps an element
$ 1 \otimes 1 \otimes r $ to $ x \otimes 1 \otimes r + 1 \otimes x \otimes r. $
On the other hand, $ {y}_i(1 \otimes 1 \otimes r) = -x \otimes 1\otimes r. $
Also, $ {y}_{i+1}(1 \otimes r) = 1 \otimes x \otimes r. $ Thus both sides are
the same
\end{itemize}

Case 6: $ (\lambda_i, \lambda_{i+1}, \lambda_{i+2})=(1,1,1). $ In this case, $
B \cong \mathfrak{F}(D_{\lambda+\alpha_{i+1},i} \circ T_{\lambda,i+1}) $ and $ B' \cong
\mathfrak{F}(D^{\lambda+\alpha_i,i+1} \circ T_{\lambda,i}). $
 Let $ a $ and $ b $ be crossingless matches.
\begin{itemize}
\item Suppose that the
circle passing through point $ i+1 $ on the bottom row of $
D_{\lambda+\alpha_{i+1},i} \circ T_{\lambda,i+1} $ is the same as the circle
passing through point $ i $ on the bottom row. Then $ {}_a B_{b} = \mathcal{A}
\otimes R $ and $ {}_a B_{b}' = \mathcal{A} \otimes \mathcal{A} \otimes R. $
Then the map on the left side of the proposition is $ (m \otimes 1)
\circ (\Delta \otimes 1). $ Thus it maps an element $ 1 \otimes r $ to $
2x \otimes r. $ On the other hand, $ {y}_i(1 \otimes r) = x \otimes r. $ Also,
$ {y}_{i+1}(1 \otimes r) = x \otimes r. $ Thus both sides are the same.
\item Suppose that the circle passing through point $ i+1 $ on the bottom row of $
D_{\lambda+\alpha_{i+1},i} \circ T_{\lambda,i+1} $ is different from the circle
passing through point $ i $ on the bottom row. Then $ {}_a B_{b} = \mathcal{A}
\otimes \mathcal{A} \otimes R $ and $ {}_a B_{b}' = \mathcal{A} \otimes R. $
Then the map on the left side of the proposition is $ (\Delta \otimes 1)
\circ (m \otimes 1). $ Thus it maps an element $ 1 \otimes 1 \otimes r $
to $ x \otimes 1 \otimes r + 1 \otimes x \otimes r. $ On the other hand, $
{y}_i(1 \otimes 1 \otimes r) = x \otimes 1\otimes r. $ Also, $ {y}_{i+1}(1
\otimes r) = 1 \otimes x \otimes r. $  Thus both sides are the same.
\end{itemize}

Case 7: $ (\lambda_i, \lambda_{i+1}, \lambda_{i+2})=(1,1,2). $
This is similar to case 5.

Case 8: $ (\lambda_i, \lambda_{i+1}, \lambda_{i+2})=(0,1,2). $
This is similar to case 6.
\end{proof}

\begin{prop}
\label{kl3.15} Let $i,j\in I^+$. If $i\neq j$, then
\begin{enumerate}
\item $ ({1}_j  {y}_i) \circ {\psi}_{i,j} = {\psi}_{i,j} \circ ({y}_i  {1}_j). $
\item $ ({y}_j  {1}_i) \circ {\psi}_{i,j} = {\psi}_{i,j} \circ ({1}_i  {y}_{j}). $
\end{enumerate}
\end{prop}

\begin{proof}
We prove only the first statement. Assume further that $ j=i+1$, the case $ j=i-1 $ being similar.  The case for $ |j-i|>1 $ is easy because the bimodules for $ \mathbb{E}_i \mathbb{E}_j $ and
$ \mathbb{E}_j \mathbb{E}_i $ are equal.

There are four non-trivial case for $ (\lambda_i, \lambda_{i+1}, \lambda_{i+2}). $ Let $ a $ and
$ b $ be crossingless matches. Let $ B $ be the bimodule for $ \mathbb{E}_i
\mathbb{E}_{i+1} $ and let $ B' $ be the bimodule for $ \mathbb{E}_{i+1}
\mathbb{E}_i. $

Case 1: $ (\lambda_i, \lambda_{i+1}, \lambda_{i+2}) = (1,1,2). $
\begin{itemize}
\item Suppose the circle passing through point $ i $ point on the
bottom row of the tangle for $ \mathbb{E}_i \mathbb{E}_{i+1} $ is the same as
the circle passing through point $ i+1 $ on the bottom row. Then $ {}_a B_b =
\mathcal{A} \otimes R $ and $ {}_a B_{b}' =  \mathcal{A} \otimes  \mathcal{A}
\otimes R $ where $ R $ denotes a tensor product of $ \mathcal{A} $
corresponding to the remaining circles. Then $ {\psi}_{i,i+1} $ is given by $
\Delta \otimes 1. $ Then $ {\psi}_{i,i+1} {y}_i(1 \otimes
r)={\psi}_{i,i+1} (x \otimes r)= x \otimes x \otimes r. $ Then $ {y}_i
{\psi}_{i,i+1} (1 \otimes r)= {y}_i (x \otimes 1 \otimes r + 1 \otimes x
\otimes r)=x \otimes x \otimes r. $
\item Suppose the circle passing through point $ i $ on the bottom row of the
tangle for $ \mathbb{E}_i \mathbb{E}_{i+1} $ is different from the circle
passing through point $ i+1 $ on the bottom row. Then $ {}_a B_b= \mathcal{A}
\otimes \mathcal{A} \otimes R $ and $ {}_a B_b'=  \mathcal{A} \otimes R. $ Then
$ {\psi}_{i,i+1} = m \otimes 1. $ Then it is easy to verify that $
{\psi}_{i,i+1} {y}_i (1 \otimes 1 \otimes r)={y}_i {\psi}_{i,i+1} (1 \otimes 1
\otimes r)=x \otimes r. $
\end{itemize}

Case 2: $ (\lambda_i, \lambda_{i+1}, \lambda_{i+2}) = (0,1,1). $
Similar to case 1.

Case 3: $ (\lambda_i, \lambda_{i+1}, \lambda_{i+2}) = (1,1,1). $
\begin{itemize}
\item Suppose the circle passing through point $ i $ on the
bottom row of the tangle is the same as the circle passing through point $ i+1
$ on the bottom row. Then $ {}_a B_b = \mathcal{A} \otimes R $ and $ {}_a
B_{b}' =  \mathcal{A} \otimes \mathcal{A} \otimes R. $ Then $ {\psi}_{i,i+1} $
is given by $ \Delta \otimes 1. $ This then follows as in case 1.
\item Suppose the circle passing through point $ i $ on the bottom row of the
tangle is different from the circle passing through the point $ i+1 $ on the
bottom row. Then $ {}_a B_b= \mathcal{A} \otimes \mathcal{A} \otimes R $ and $
{}_a B_b'=  \mathcal{A} \otimes R. $ Then $ {\psi}_{i,i+1} = m \otimes
1. $ This then follows as in case 1.
\end{itemize}

Case 4: $ (\lambda_i, \lambda_{i+1}, \lambda_{i+2}) = (0,1,2). $
Similar to case 3.
\end{proof}

\begin{prop}\label{kl3.17} For $i,j,k\in I^+$,
$$(\psi_{j,k}1_i)\circ (1_j\psi_{i,k})\circ(\psi_{i,j}1_k)+(1_k\psi_{i,j})\circ(\psi_{i,k}1_j)\circ(1_i \psi_{j,k})
=\begin{cases}0&i\neq k\mbox{ or }|i-j| \neq 1,\\  1_i  1_j  1_i&i=k\mbox{ and
}|i-j|=1.\end{cases}$$

\end{prop}

\begin{proof}
The proof of the first part consists of verifying the equality in many different cases, each of which is similar to the second part.
We only prove the second part in the case $ j = i+1 $ as the case $ j = i-1 $ is similar.  There are four cases  for $ (\lambda_i, \lambda_{i+1}, \lambda_{i+2}) $ for which $ \mathbb{E}_i \mathbb{E}_{i+1} \mathbb{E}_i $ is non-zero.

Case 1:  $ (\lambda_i, \lambda_{i+1}, \lambda_{i+2})=(0,1,1). $  In this case,
$ ({\psi}_{j,i}  {1}_i) \circ ({1}_j  {\psi}_{i,i}) \circ ({\psi}_{i,j}
 {1}_i)  = 0 $ because it passes through the functor $ \mathbb{E}_{i+1}
\mathbb{E}_i \mathbb{E}_i $ which is zero on the category corresponding to this
$ \lambda. $ On the other hand
$$ ({1}_i  {\psi}_{i,j}) \circ ({\psi}_{i,i}  {1}_j) \circ ({1}_i  {\psi}_{j,i}) \colon
\mathbb{E}_i \mathbb{E}_{i+1} \mathbb{E}_i \rightarrow \mathbb{E}_i \mathbb{E}_{i} \mathbb{E}_{i+1} \rightarrow \mathbb{E}_i \mathbb{E}_{i} \mathbb{E}_{i+1} \rightarrow \mathbb{E}_i \mathbb{E}_{i+1} \mathbb{E}_i. $$
Let $ B $ be the bimodule for the functor $ \mathbb{E}_i \mathbb{E}_{i+1} \mathbb{E}_i. $
Then this is a sequence of maps
$$ B \rightarrow \mathcal{A} \otimes B \rightarrow \mathcal{A} \otimes B \rightarrow B $$
where the first map given by comultiplication, the middle map is given by the map $ 1\otimes\kappa,$ and the last map is multiplication.
This sequence of maps acts on $ 1 \otimes \alpha \in B $ as follows:
$$ 1 \otimes \alpha \mapsto x \otimes 1 \otimes \alpha + 1 \otimes x \otimes \alpha \mapsto 1 \otimes 1 \otimes \alpha \mapsto 1 \otimes \alpha. $$
Clearly $ (({\psi}_{j,i}  {1}_i) \circ ({1}_j  {\psi}_{i,i}) \circ ({\psi}_{i,j}  {1}_i) (1 \otimes \alpha) =0. $
Similarly, $ x \otimes \alpha \mapsto x \otimes \alpha. $

Case 2: $ (\lambda_i, \lambda_{i+1}, \lambda_{i+2})=(0,2,2). $ This is similar
to case 1 except that now $({1}_i  {\psi}_{i,j}) \circ ({\psi}_{i,i} {1}_j)
\circ ({1}_i  {\psi}_{j,i})=0$ and $ ({\psi}_{j,i}  {1}_i) \circ ({1}_j
 {\psi}_{i,i}) \circ ({\psi}_{i,j}  {1}_i) = 1_i 1_j 1_i. $

Case 3: $ (\lambda_i, \lambda_{i+1}, \lambda_{i+2})=(0,1,2). $  In this case, $
({\psi}_{j,i}  {1}_i) \circ ({1}_j  {\psi}_{i,i}) \circ ({\psi}_{i,j}
 {1}_i)  = 0 $ since this map passes through the functor $
\mathbb{E}_{i+1} \mathbb{E}_i \mathbb{E}_i $ which is zero on the category
corresponding to $ \lambda. $

On the other hand
$$ ({1}_i  {\psi}_{i,j}) \circ ({\psi}_{i,i}  {1}_j) \circ ({1}_i  {\psi}_{j,i}) \colon
\mathbb{E}_i \mathbb{E}_{i+1} \mathbb{E}_i \rightarrow \mathbb{E}_i \mathbb{E}_{i} \mathbb{E}_{i+1} \rightarrow \mathbb{E}_i \mathbb{E}_{i} \mathbb{E}_{i+1} \rightarrow \mathbb{E}_i \mathbb{E}_{i+1} \mathbb{E}_i. $$
Let $ B $ be the bimodule for the functor $ \mathbb{E}_i \mathbb{E}_{i+1} \mathbb{E}_i. $
Then this is a sequence of maps
$$ B \rightarrow \mathcal{A} \otimes B \rightarrow \mathcal{A} \otimes B \rightarrow B $$
where the first and third maps are given by lemmas ~\ref{mulemma1} and ~\ref{mulemma2} respectively, and the middle map is given in section ~\ref{2morphisms}.
This sequence of maps acts on $ 1 \otimes \alpha, x\otimes\alpha \in B $ as follows:
 \begin{align}1 \otimes \alpha \mapsto x \otimes 1 \otimes \alpha + 1 \otimes x \otimes \alpha \mapsto 1 \otimes 1 \otimes  \alpha \mapsto 1 \otimes \alpha,
 x\otimes\af\mapsto x\otimes x\otimes\af\mapsto x\otimes 1\otimes \af\mapsto x\otimes\af. \end{align}

Case 4: $ (\lambda_i, \lambda_{i+1}, \lambda_{i+2})=(0,2,1). $ This is similar
to case 1 except that now $ ({1}_i  {\psi}_{i,j}) \circ ({\psi}_{i,i} {1}_j)
\circ ({1}_i  {\psi}_{j,i}) = 0 $ and $ ({\psi}_{j,i}  {1}_i) \circ ({1}_j
 {\psi}_{i,i}) \circ ({\psi}_{i,j}  {1}_i) (\beta \otimes \alpha) = \beta \otimes \alpha. $
\end{proof}

\begin{theorem} There is a 2-functor $ \Omega_{k,n} \colon \mathcal{KL} \rightarrow \mathcal{HK}_{k,n} $ such that for all $i,j\in I$,
\begin{enumerate}
\item $ \Omega_{k,n}(\lambda) = \mathcal{C}_{\lambda}, $
\item $ \Omega_{k,n}(\mathcal{I}_{\ld}) = \bbI_{\ld}, $
\item $ \Omega_{k,n}(\mathcal{E}_i \mathcal{I}_{\ld}) = \mathbb{E}_i \bbI_{\ld},$
\item $ \Omega_{k,n}(Y_{i;\ld}) = {y}_{i;\ld}, $
\item $ \Omega_{k,n}(\Psi_{i,j;\ld}) = {\psi}_{i,j;\ld},$
\item $ \Omega_{k,n}(\Cup_{i;\ld}) = {\cup}_{i;\ld}, $
\item $ \Omega_{k,n}(\Cap_{i;\ld}) = {\cap}_{i;\ld}, $
\item $ \Omega_{k,n}(\bone_{i;\ld}) = 1_{i;\ld} .$
\end{enumerate}
\end{theorem}

\section{The 2-category $ \mathcal{P}_{k,n}$}\label{S:Cat O}

\subsection{Graded category $ {}_{\mathbb{Z}} \mathcal{O}$}
Let $ \mathfrak{g} = \mathfrak{gl}_{2k} $ be the Lie algebra of $2k\times2k$-matrices, let $\mathfrak{d}$ denote the Cartan subalgebra of $\g$ consisting of diagonal matrices and $ \mathfrak{p} $ be the Borel subalgebra of upper triangular matrices. For $i=1,\ldots,2k$, let $e_{ij}$ denote the $(i,j)$-matrix unit, and let $\ep_i\in \mathfrak{d}^*$ be the coordinate functional $\ep_i(e_{jj})=\dt_{ij}$. Let $\mathcal{O} $ be the category of finitely generated
$ \mathfrak{g}$-modules which are diagonalizable with respect to $\mathfrak{d}$ and locally finite with respect to $ \mathfrak{p} $. Let
\[
X=\bigoplus_{i=1}^{2k}\Z\ep_i,\;\;\;\mbox{and}\;\;\;Y=\bigoplus_{i=1}^{{2k}-1}\Z(\ep_i-\ep_{i+1})\subset
X
\]
denote the weight lattice and root lattice of $\gl_{{2k}}$, respectively. The
dominant weights are given by the set
$X^+=\{\,\mu=\mu_1\ep_1+\cdots+\mu_{2k}\ep_{2k}\in
X\,|\,\mu_1\geq\cdots\geq\mu_{2k}\,\}$. Denote half the sum of the positive roots by $ \rho. $
Let $ \mu\in X^+ $, and $\mathcal{O}_{\mu}$ the block of $\mathcal{O}$ consisting of modules that have a generalized central character corresponding to
$\mu$ under the Harish-Chandra homomorphism.  Let $\mathcal{O}^{(k,k)}_{\mu}$ be the full subcategory $\mathcal{O}$ consisting of modules which are locally finite with respect to the parabolic subalgebra whose reductive part is
$ \mathfrak{gl}_k \oplus \mathfrak{gl}_k. $
Finally, let $ \mathcal{P}_{\mu}^{(k,k)} $ be the full subcategory of $ \mathcal{O}^{(k,k)}_{\mu} $ whose objects have projective presentations by projective-injective modules.

Let $\mu$ and $\mu'$ be integral dominant weights of $\g$, and let $ \Stab(\mu) $ denote the stabilizer of $ \mu $ under the $\rho$-shifted action of the symmetric group $ \mathbb{S}_{2k}$.
Suppose $ \mu'-\mu $ is an integral dominant weight.
Then, let $ \theta_{\mu}^{\mu'} \colon \mathcal{O}^{(k,k)}_{\mu} \rightarrow \mathcal{O}^{(k,k)}_{\mu'} $ be the translation functor of tensoring with the finite dimensional irreducible representation of highest weight $ \mu' - \mu $ composed with projecting onto the $\mu'$-block,
and let $ \theta_{\mu'}^{{\mu}} $ be its adjoint.

Let $ P_{\mu} $ be a minimal projective generator of $ \mathcal{O}_{\mu}. $ It
was shown that $ A_{\mu} = \End_{\mathfrak{g}}(P_{\mu}) $ has the structure of
a graded algebra \cite{bgs}. Since $ \mathcal{O}_{\mu} $ is Morita equivalent
to $ A_{\mu}$-mod, we consider the category of graded $ A_{\mu}$-modules which
we denote by $ {}_{\mathbb{Z}} \mathcal{O}_{\mu}. $ Let the graded lift of $
\mathcal{O}^{(k,k)}_{\mu} $ and $ \mathcal{P}^{(k,k)}_{\mu} $ be $
{}_{\mathbb{Z}} \mathcal{O}^{(k,k)}_{\mu} $ and $ {}_{\mathbb{Z}}
\mathcal{P}^{(k,k)}_{\mu} $, respectively. It is known that if $ \Stab(\mu)
\subset \Stab(\mu'),$ there is a graded lift of the translation functors, cf.
\cite{str1}, which by abuse of notation we denote again by $
\widetilde{\theta}_{\mu'}^{\mu} $ and $ \widetilde{\theta}_{\mu}^{\mu'} $.

The key tool in the construction of graded category $\mathcal{O}$ is the
Soergel functor. Let $\ld=(\ld_1,\ldots,\ld_n)$ be a composition of $2k$, let
$\mathbb{S}_{\ld} = \mathbb{S}_{\ld_1} \times \cdots \times
\mathbb{S}_{\ld_{n}},$ let $w^\mu_0$ be the longest coset representative in
$\mathbb{S}_{2k}/\mathbb{S}_\mu$, and let $P(w^\mu_0\cdot\mu)$ be the unique up
to isomorphism, indecomposable projective-injective object of
$\mathcal{O}_{\mu}. $ Let
$C=S(\mathfrak{h})/S(\mathfrak{h})_+^{\mathbb{S}_{2k}} $ be the coinvariant
algebra of the symmetric algebra for the Cartan subalgebra with respect to the
action of the symmetric group. Let $ x_1, \ldots, x_{2k} $ be a basis of
$S(\mathfrak{h})$ and by abuse of notation also let $ x_i $ denote its image in
$ C. $ Let $ C^{\ld} $ be the subalgebra of elements invariant under the action
of $\mathbb{S}_{\ld}$.  Soergel proved in \cite{soe1}:

\begin{prop}
$ \End_{\mathfrak{g}}(P(w^\mu_0\cdot\mu)) \cong C^{\Stab(\mu)}. $
\end{prop}

Define the Soergel functor $ \mathbb{V}_{\mu} \colon \mathcal{O}_{\mu}
\rightarrow C^{\Stab(\mu)}$-mod to be $ \Hom_{\mathfrak{g}}(P(w_0.\mu),
\bullet). $

\begin{prop}
\label{fullyfaithful} Let $ P $ be a projective object.  Then there is a natural isomorphism $
\Hom_{C^{\Stab(\mu)}}(\mathbb{V}_{\mu}P, \mathbb{V}_{\mu} M) \cong
\Hom_{\mathfrak{g}}(P, M). $
\end{prop}

\begin{proof}
This is the Structure Theorem of \cite{soe1}.
\end{proof}

\begin{prop}\label{isom of functors}
Let $ \mu,\mu'\in X^+ $ be integral dominant weights such that there is a
containment of stabilizers: $ \Stab({\mu}) \subset \Stab({\mu'}). $ Then there
are isomorphism of functors
\begin{enumerate}
\item $ \mathbb{V}_{\mu'} \theta_{\mu}^{\mu'} \cong \Res_{C^{\Stab(\mu)}}^{C^{\Stab(\mu')}} \mathbb{V}_{\mu} $
\item $ \mathbb{V}_{\mu} \theta_{\mu'}^{\mu} \cong C^{\Stab(\mu)} \otimes_{C^{\Stab(\mu')}} \mathbb{V}_{\mu'}. $
\end{enumerate}
\end{prop}

\begin{proof}
This is Theorem 12 and Proposition 6 of \cite{soe2}.
\end{proof}

\subsection{The objects of $ \mathcal{P}_{k,n}$}
Let $ \lambda = (\lambda_1, \ldots, \lambda_r) $ be a composition of $ 2k $ with $ \ld_i \in \lbrace 0,1,2 \rbrace $ for all $ i. $
To each such $ \lambda, $ we associate an integral dominant weight
$$
\displaystyle
\overline{\lambda} = \sum_{j=1}^r \sum_{i=1}^{\lambda_j} (r-j+1)\ep_{\ld_1+\cdots+\lambda_{j-1}+i} - \rho $$
of $ \mathfrak{gl}_{2k} $ where $ \lambda_0=0. $ Note the stabilizer of this weight under the action of $ \mathbb{S}_{2k} $ is
$ \mathbb{S}_{\lambda_1} \times \cdots \times \mathbb{S}_{\lambda_n}. $

The set of objects of $ \mathcal{P}_{k,n} $ are the categories $ {}_{\mathbb{Z}} \mathcal{P}_{\overline{\lambda}}^{(k,k)}, $ $\ld\in P(V_{2\om_k})$.

\subsection{The 1-morphisms of $ \mathcal{P}_{k,n}$}
Let $\ld\in P(V_{2\om_k})$, and let $ \mathbb{I}_{\ld}\in\End_\g({}_{\mathbb{Z}} \mathcal{P}_{\overline{\lambda}}^{(k,k)})$ be the identity functor.

For each $i\in I$, we define functors $\mathbb{E}_i \mathbb{I}_{\ld}$, and $\bbK_i\bbI_\ld$
To this end, let $ \lambda $ be a weight of $ V_{2 \omega_k} $ and $ i\in I^+. $  Then we have compositions of $ 2k $ into $ n+1 $ parts:
\[
\ld(i)=(\ld_1, \ldots, \ld_i, 1, \ld_{i+1}-1, \ldots, \ld_n)      ,\;\;\; \ld(-i)=(\ld_1, \ldots, \ld_i-1, 1, \ld_{i+1}, \ldots, \ld_n)
\]
Also, if $\ld=\sum_ia_i\om_i\in P$, set $r_{i,\ld}=1+a_1+\cdots+a_{i-1}+a_{i+1} $ and $ s_{i,\ld}=2-a_i-a_{i+1}. $

Let $i\in I$. Suppose $ (\lambda_i, \lambda_{i+1}) \in \lbrace (0,1), (0,2), (1,1), (1,2) \rbrace. $  Then we define as in \cite{fks}, $ \mathbb{E}_i\mathbb{I}_{\ld} \colon {}_{\mathbb{Z}} \mathcal{P}_{\overline{\lambda}}^{(k,k)} \rightarrow {}_{\mathbb{Z}}\mathcal{P}_{\overline{\lambda + \alpha_i}}^{(k,k)} $ is given by tensoring with the following bimodule:
\begin{align*}
\Hom_{\mathfrak{g}}(P_{\overline{\lambda + \alpha_i}}, {\theta}_{\overline{\lambda(i)}}^{\overline{\lambda+ \alpha_i}} {\theta}_{\overline{\lambda}}^{\overline{\lambda(i)}}P_{\overline{\lambda}} \lbrace r_{i,\ld} \rbrace)
&\cong \Hom_{C^{{\lambda+\alpha_i}}}(\mathbb{V}_{\overline{\lambda + \alpha_i}} P_{\overline{\lambda + \alpha_i}}, \mathbb{V}_{\overline{\lambda + \alpha_i}} {\theta}_{\overline{\lambda(i)}}^{\overline{\lambda+ \alpha_i}} {\theta}_{\overline{\lambda}}^{\overline{\lambda(i)}}P_{\overline{\lambda}}\lbrace r_{i,\ld} \rbrace )\\
&\cong \Hom_{C^{{\lambda+\alpha_i}}}(\mathbb{V}_{\overline{\lambda + \alpha_i}} P_{\overline{\lambda + \alpha_i}}, C^{{\lambda+\alpha_i}} \otimes_{C^{{\lambda(i)}}} \Res_{C^{{\lambda}}}^{C^{{\lambda(i)}}} \mathbb{V}_{\overline{\lambda}}P_{\overline{\lambda}}\lbrace r_{i,\ld} \rbrace).
\end{align*}
For all other values of $ (\lambda_i, \lambda_{i+1}), $ set $ \mathbb{E}_i \mathbb{I}_{\ld}= 0. $
Let $ \mathbb{K}_i\mathbb{I}_{\ld} \colon {}_{\mathbb{Z}} \mathcal{P}_{\overline{\lambda}}^{(k,k)} \rightarrow {}_{\mathbb{Z}} \mathcal{P}_{\overline{\lambda}}^{(k,k)} $ be the grading shift functor $\bbK_i\bbI_\ld=\bbI_\ld\{(\af_i,\ld)\}$.

Let $ {}_{\mathbb{Z}}\mathcal{P}_{\overline{\lambda}}^{(k,k)} $ and $ {}_{\mathbb{Z}}\mathcal{P}_{\overline{\lambda}'}^{(k,k)} $
be two objects.  Then
\[
\Hom({}_{\mathbb{Z}}\mathcal{P}_{\overline{\lambda}}^{(k,k)}, {}_{\mathbb{Z}}\mathcal{P}_{\overline{\lambda}'}^{(k,k)})
=\bigoplus_{\substack{\ui\in\Seq\\s\in\Z}}\C
\bbI_{\ld'}\bbE_{\ui}\bbI_{\ld}\{s\}
\]
where $\bbE_\ui:=\bbE_{i_1}\cdots \bbE_{i_r}\bbI_{\ld}$ if $\ui=(i_1,\ldots,i_r)\in I_\infty$,
and $s$ refers to a \emph{grading shift}.

\subsection{Bimodule categories over the cohomology of flag varieties}
A review of certain bimodules and bimodule maps over the cohomology of flag varieties developed in \cite{l, kl, b} is given here.
Let $ \ld = (\ld_1, \ldots, \ld_n) $ be a composition of $ 2k $ into $ n $ parts.  Let $ x(\ld)_{j,r}=x_{\ld_1+\cdots+\ld_{j-1}+r}. $
There is an isomorphism of algebras:
$$ \displaystyle
C^{{\ld}} \cong \bigotimes_{1\leq j\leq n} \mathbb{C}[x(\ld)_{j,1}, x(\ld)_{j,2}, \ldots x(\ld)_{j,\ld_j}] / J_{\ld, {n}} $$
where $ J_{\ld, {n}} $ is the ideal generated by the homogeneous terms in the equation
\begin{equation}\label{E:Homogeneous Terms}
\prod_{1\leq j\leq n} (1+ x(\ld)_{j,1} t + x(\ld)_{j,2} t^2 + \cdots + x(\ld)_{j,\ld_j} t^{\ld_j}) = 1.
\end{equation}

Let $\overline{x(\ld)_{i,k}}$ be the homogenous term of degree $2k$ in the product $$\prod_{\substack{1\leq j\leq n\\j\neq i}} (1+ x(\ld)_{j,1} t + x(\ld)_{j,2} t^2 + \cdots + x(\ld)_{j,\ld_j} t^{\ld_j}).$$ Then, using \eqref{E:Homogeneous Terms} we see that $$\sum_{j=1}^{k}x(\ld)_{i,j}\overline{x(\ld)_{i,k-j}}=\dt_{k,0},$$ cf. \cite[$\S5.1$]{kl} for details.

We must also consider $ C^{{\ld(i)}}. $
There is an isomorphism of algebras:
$$ C^{{\ld(i)}} \cong \bigotimes_{\substack{1\leq j\leq n,\\j\neq i+1}} \mathbb{C}[x(\ld)_{j,1}, x(\ld)_{j,2}, \ldots, x(\ld)_{j,\ld_j}] \otimes \mathbb{C}[\zeta_i] \otimes
\mathbb{C}[x(\ld)_{i+1,1}, x(\ld)_{i+1,2}, \ldots, x(\ld)_{i+1,\ld_{i+1}-1}]
/ J_{\ld(i), {n}} $$
where $ J_{\ld(i), {n}} $ is the ideal generated by the homogeneous terms in the equation
$$ \displaystyle
\prod_{\substack{1\leq j\leq n, \\j \neq i+1}} (1+\zeta_i t)
\sum_{r=0}^{\ld_{i+1}-1} x(\ld)_{i+1,r}t^r
\sum_{s=0}^{\ld_j} x(\ld)_{j,s}t^s
= 1 .$$
There is also an isomorphism of algebras:
$$ C^{{\ld(-i)}} \cong \bigotimes_{\substack{1\leq j\leq n,\\j\neq i}} \mathbb{C}[x(\ld)_{j,1}, x(\ld)_{j,2}, \ldots, x(\ld)_{j,\ld_j}] \otimes
\mathbb{C}[x(\ld)_{i,1}, x(\ld)_{i,2}, \ldots, x(\ld)_{i,\ld_{i}-1}] \otimes \mathbb{C}[\zeta_i]
/ J_{\ld(-i), {n}} $$
where $ J_{\ld(-i), {n}} $ is the ideal generated by the homogeneous terms in the equation
$$ \displaystyle
\prod_{\substack{1\leq j\leq n, \\j \neq i}} (1+\zeta_i t)
\sum_{r=0}^{\ld_{i}-1} x(\ld)_{i,r}t^r
\sum_{s=0}^{\ld_j} x(\ld)_{j,s}t^s
= 1 .$$

\subsection{The $2$-morphisms} In light of Propositions \ref{fullyfaithful} and \ref{isom of functors}, we may define the 2-morphisms on the algebras $C^{{\ld}}$, $\ld\in P(V_{2\om_k})$.

\subsection*{The maps $ {\overline{y}_{i;\ld}}$}
Let $i\in I$. Define  $ \overline{y}_{i;\ld} \colon
C^{{\lambda(i)}} \rightarrow C^{{\lambda(i)}} $ which is a map of $
(C^{{\lambda+\alpha_i}}, C^{{\lambda}}) $-bimodules by $ \overline{y}_{i;\ld}((\zeta_i)^r) =
(\zeta_{i})^{r+1}. $

\subsection*{The maps $ \overline{\cup}_{i;\ld}, \overline{\cap}_{i;\ld}$} Let $i\in I^+$. Define a map of $ (C^{{\lambda}}, C^{{\lambda}})$-bimodules
$$ \overline{\cup}_{i;\ld} \colon C^{{\lambda}} \rightarrow C^{{\lambda(i)}} \otimes_{C^{{\lambda + \alpha_i}}} C^{{\lambda(i)}} \lbrace 1-\ld_i-\ld_{i+1} \rbrace $$
by
$$ \displaystyle
\overline{\cup}_{i;\ld}(1)=\sum_{f=0}^{\lambda_{i}}
(-1)^{\lambda_{i} -f} \zeta_i^f \otimes
x(\lambda)_{i,\lambda_{i}-f}. $$

Next define a map of $ (C^{{\lambda}}, C^{{\lambda}})$-bimodules
$$ \overline{\cup}_{-i;\ld} \colon C^{{\lambda}} \rightarrow C^{{\lambda(-i)}} \otimes_{C^{{\lambda - \alpha_i}}} C^{{\lambda(-i)}} \lbrace 1-\ld_i-\ld_{i+1} \rbrace $$
by
$$ \displaystyle
\overline{\cup}_{-i;\ld}(1)=\sum_{f=0}^{\lambda_{i+1}}
(-1)^{\lambda_{i+1}-f} \zeta_i^f \otimes
x(\lambda)_{i+1,\lambda_{i+1}-f}. $$

Next define a map of $ (C^{{\lambda}}, C^{{\lambda}})$-bimodules
$$ \overline{\cap}_{i;\ld} \colon C^{{\lambda(i)}} \otimes_{C^{{\lambda + \alpha_i}}} C^{{\lambda(i)}} \lbrace 1-\ld_i-\ld_{i+1} \rbrace \rightarrow C^{{\lambda}} $$
by
$$ \overline{\cap}_{i;\ld}(\zeta_i^{r_1} \otimes \zeta_i^{r_2})=
(-1)^{r_1+r_2+1-\lambda_{i+1}} \overline{x(\lambda)_{i+1,r_1+r_2+1-\lambda_{i+1}}}. $$

Next define a map of $ (C^{{\lambda}}, C^{{\lambda}})$-bimodules
$$ \overline{\cap}_{-i;\ld} \colon C^{{\lambda(-i)}} \otimes_{C^{{\lambda - \alpha_i}}} C^{{\lambda(-i)}} \lbrace 1-\ld_i-\ld_{i+1} \rbrace \rightarrow C^{{\lambda}} $$
by
$$ \overline{\cap}_{-i;\ld}(\zeta_i^{r_1} \otimes \zeta_i^{r_2})=
(-1)^{r_1+r_2+1-\lambda_{i}} \overline{x(\lambda)_{i,r_1+r_2+1-\lambda_{i}}}. $$

\subsection*{The maps $ \overline{\psi}_{i,j;\ld} $} Let $i,j\in I^+$.
Define a map of $ (C^{{\lambda+\alpha_i+\alpha_j}}, C^{{\lambda}})$-bimodules
$$ \displaystyle
\overline{\psi}_{i,j;\ld} \colon C^{{(\lambda+\alpha_j)(i)}}
\otimes_{{C^{\lambda+\alpha_j}}} C^{{\lambda(j)}} \rightarrow
C^{{(\lambda+\alpha_i)(j)}}  \otimes_{C^{{\lambda+\alpha_i}}} C^{{\lambda(i)}} $$
by
\begin{eqnarray*}
\overline{\psi}_{i,j;\ld}(\zeta_i^{r_1} \otimes \zeta_j^{r_2}) =\begin{cases}
\zeta_j^{r_2} \otimes \zeta_i^{r_1} & \textrm{if }  |i-j|>1, \\
\sum_{f=0}^{r_1-1} \zeta_i^{r_1+r_2-1-f} \otimes \zeta_i^{f} - \sum_{g=0}^{r_2-1} \zeta_i^{r_1+r_2-1-g} \otimes \zeta_i^{g} & \textrm{if } j=i, \\
(\zeta_j^{r_2} \otimes \zeta_i^{r_1+1} - \zeta_j^{r_2+1} \otimes \zeta_i^{r_1}) \lbrace -1 \rbrace & \textrm{if } i = j+1, \\
(\zeta_j^{r_2} \otimes \zeta_i^{r_1}) \lbrace 1 \rbrace & \textrm{if } j = i+1.\end{cases}
\end{eqnarray*}

Define a map of $ (C^{{\lambda-\alpha_i-\alpha_j}}, C^{{\lambda}})$-bimodules
$$ \overline{\psi}_{-i,-j;\ld} \colon C^{{(\lambda-\alpha_j)(-i)}}  \otimes_{C^{{\lambda-\alpha_j}}} C^{{\lambda(-j)}}
\rightarrow C^{{(\lambda-\alpha_i)(-j)}}  \otimes_{C^{{\lambda-\alpha_i}}}
C^{{\lambda(-i)}} $$ by

\begin{eqnarray*}
\overline{\psi}_{-i,-j}(\zeta_i^{r_1} \otimes \zeta_j^{r_2}) =\begin{cases}
\zeta_j^{r_2} \otimes \zeta_i^{r_1} & \textrm{if }  |i-j|>1,\\
\sum_{f=0}^{r_2-1} \zeta_i^{r_1+r_2-1-f} \otimes \zeta_i^{f} - \sum_{g=0}^{r_1-1} \zeta_i^{r_1+r_2-1-g} \otimes \zeta_i^{g} & \textrm{if } j=i, \\
(\zeta_j^{r_2} \otimes \zeta_i^{r_1+1}) \lbrace -1 \rbrace  & \textrm{if } i = j+1, \\
(\zeta_j^{r_2+1} \otimes \zeta_i^{r_1}- \zeta_j^{r_2} \otimes \zeta_i^{r_1+1}) \lbrace 1 \rbrace & \textrm{if } j = i+1.\end{cases}
\end{eqnarray*}

\subsection{The 2-morphisms of $ \mathcal{P}_{k,n}$}
Let $ i,j \in I^{+}. $

\subsection*{The maps $ 1_{i;\ld}$}

Let $ 1_{i;\ld} \colon \mathbb{E}_i \mathbb{I}_{\ld} \rightarrow \mathbb{E}_i \mathbb{I}_{\ld}$ and $ 1_{-i;\ld} \colon
\bbE_{-i} \mathbb{I}_{\ld} \rightarrow \bbE_{-i} \mathbb{I}_{\ld} $ be the identity morphisms.

\subsection*{The maps $ y_{i;\ld}$}

Next we define a morphism of degree $ 2, $ $ y_{i;\ld} \colon \mathbb{E}_i \mathbb{I}_{\ld}
\rightarrow \mathbb{E}_i \mathbb{I}_{\ld}. $ Recall that
$$ \mathbb{E}_i \mathbb{I}_{\ld} \cong
\Hom_{C^{{\lambda+\alpha_i}}}(\mathbb{V}_{\overline{\lambda + \alpha_i}}
P_{\overline{\lambda +\alpha_i}}, C^{{\lambda(i)}} \otimes_{C^{{\lambda}}}
\mathbb{V}_{\overline{\lambda}} P_{\overline{\lambda}} \lbrace r_{i,\ld} \rbrace). $$ Let $ f $ be such a
homomorphism. Suppose $ f(m) = \gamma \otimes n. $ Then set $
(y_{i;\ld}.f)(m)=\overline{y}_i(\gamma) \otimes n . $

Similarly,
$$ \bbE_{-i} \mathbb{I}_{\ld} \cong
\Hom_{C^{{\lambda-\alpha_i}}}(\mathbb{V}_{\overline{\lambda - \alpha_i}}
P_{\overline{\lambda -\alpha_i}}, C^{{\lambda(-i)}} \otimes_{C^{{\lambda}}}
\mathbb{V}_{\overline{\lambda}}P_{\overline{\lambda}} \lbrace s_{i,\ld} \rbrace). $$ Let $ f $ be such a
homomorphism. Suppose $ f(m) = \gamma \otimes n. $ Then set $
(y_{-i;\ld}.f)(m)=\overline{y}_{-i;\ld}(\gamma) \otimes n. $

\subsection*{The maps $ \cup_{i;\ld}, \cap_{i;\ld}$}

Note that
$$ \mathbb{I}_{\ld} \cong J= \Hom_{C^{{\lambda}}}(\mathbb{V}_{\overline{\lambda}} P_{\overline{\lambda}}, \mathbb{V}_{\overline{\lambda}} P_{\overline{\lambda}}) $$
$$ \bbE_{-i} \circ \mathbb{E}_i \mathbb{I}_{\ld} \cong
K= \Hom_{C^{{\lambda}}}(\mathbb{V}_{\overline{\lambda}} P_{\overline{\lambda}}, C^{{{\lambda+\alpha_i(-i)}}} \otimes_{C^{{\lambda + \alpha_i}}} C^{{\lambda(i)}} \otimes_{C^{{\lambda}}} \mathbb{V}_{\overline{\lambda}} P_{\overline{\lambda}} \lbrace r_{\ld,i} + s_{\ld+\alpha_i,i} \rbrace) $$
$$ \mathbb{E}_i \circ \bbE_{-i} \mathbb{I}_{\ld} \cong
L=\Hom_{C^{{\lambda}}}(\mathbb{V}_{\overline{\lambda}} P_{\overline{\lambda}}, C^{{{\lambda-\alpha_i}(i)}} \otimes_{C^{{\lambda - \alpha_i}}} C^{{\lambda(-i)}} \otimes_{C^{{\lambda}}} \mathbb{V}_{\overline{\lambda}} P_{\overline{\lambda}} \lbrace s_{\ld,i} + r_{\ld-\alpha_i,i} \rbrace). $$

Let $ f \in J. $ Then define $ \cup_{i;\ld} \colon \mathbb{I}_{\ld} \rightarrow \bbE_{-i}
\mathbb{E}_i \mathbb{I}_{\ld} $ by
$$ \cup_{i;\ld}(f)(m) = \overline{\cup}_{i;\ld}(1) \otimes f(m) $$
and $ \cup_{-i;\ld} \colon \mathbb{I}_{\ld} \rightarrow \mathbb{E}_i \bbE_{-i} \mathbb{I}_{\ld} $ by
$$ \cup_{-i;\ld}(f)(m) = \overline{\cup}_{-i;\ld}(1) \otimes f(m). $$

Now define $ \cap_{i;\ld} \colon \bbE_{-i} \mathbb{E}_i \mathbb{I}_{\ld} \rightarrow \mathbb{I}_{\ld}. $ Suppose $
f \in K $ such that $ f(m)=\gamma \otimes n. $ Then set $
\cap_{i;\ld}(f)(m)=\overline{\cap}_{i;\ld}(\gamma) \otimes n. $

Next define $ \cap_{-i;\ld} \colon \mathbb{E}_i \bbE_{-i} \mathbb{I}_{\ld} \rightarrow \mathbb{I}_{\ld}. $ Suppose
$ f \in L $ such that $ f(m)=\gamma \otimes n. $ Then set $
\cap_{-i;\ld}(f)(m)=\overline{\cap}_{-i;\ld}(\gamma) \otimes n. $

\subsection*{The maps $ \psi_{i,j;\ld}^{}$}

First we define a map $ \psi_{i,j;\ld} \colon \mathbb{E}_i \mathbb{E}_j \mathbb{I}_{\ld} \rightarrow
\mathbb{E}_j \mathbb{E}_i \mathbb{I}_{\ld}. $

Set
$$ J_{i,j}^+ = \mathbb{E}_i \mathbb{E}_j \mathbb{I}_{\ld} \cong
\Hom_{C^{{\lambda+\alpha_i+\alpha_j}}}(\mathbb{V}_{\overline{\lambda+\alpha_i+\alpha_j}} P_{\overline{\lambda+\alpha_i+\alpha_j}},
C^{{(\lambda+\alpha_j)(i)}}  \otimes_{C^{{\lambda+\alpha_j}}} C^{{\lambda(j)}} \otimes_{C^{{\lambda}}} \mathbb{V}_{\overline{\lambda}} P_{\overline{\lambda}} \lbrace r_{\ld,j}+r_{\ld+\alpha_j,i} \rbrace) $$
$$ K_{i,j}^+ = \mathbb{E}_j \mathbb{E}_i \mathbb{I}_{\ld} \cong
\Hom_{C^{{\lambda+\alpha_j+\alpha_i}}}(\mathbb{V}_{\overline{\lambda+\alpha_j+\alpha_i}} P_{\overline{\lambda+\alpha_j+\alpha_i}},
C^{{(\lambda+\alpha_i)(j)}}  \otimes_{C^{{\lambda+\alpha_i}}} C^{{\lambda(i)}} \otimes_{C^{{\lambda}}} \mathbb{V}_{\overline{\lambda}} P_{\overline{\lambda}} \lbrace r_{\ld,i}+r_{\ld+\alpha_i,j} \rbrace). $$

Let $ f \in J_{i,j}^+$  and suppose that $ f(m)= \gamma_1 \otimes \gamma_2
\otimes n. $ Then define $ \psi_{i,j;\ld} f(m) = \overline{\psi}_{i,j;\ld}(\gamma_1 \otimes
\gamma_2) \otimes n. $

Set
$$ J_{i,j}^- = \bbE_{-i} \mathbb{E}_{-j} \mathbb{I}_{\ld} \cong
\Hom_{C^{{\lambda-\alpha_i-\alpha_j}}}(\mathbb{V}_{\overline{\lambda-\alpha_i-\alpha_j}} P_{\overline{\lambda-\alpha_i-\alpha_j}},
C^{{(\lambda-\alpha_j)(-i)}}  \otimes_{C^{{\lambda-\alpha_j}}} C^{{\lambda(-j)}} \otimes_{C^{{\lambda}}} \mathbb{V}_{\overline{\lambda}} P_{\overline{\lambda}} \lbrace s_{\ld,j} +s_{\ld-\alpha_j,i} \rbrace) $$
$$ K_{i,j}^- = \mathbb{E}_{-j} \bbE_{-i} \mathbb{I}_{\ld} \cong
\Hom_{C^{{\lambda-\alpha_j-\alpha_i}}}(\mathbb{V}_{\overline{\lambda-\alpha_j-\alpha_i}} P_{\overline{\lambda-\alpha_j-\alpha_i}},
C^{{(\lambda-\alpha_i)(-j)}}  \otimes_{C^{{\lambda-\alpha_i}}} C^{{\lambda(-i)}} \otimes_{C^{{\lambda}}} \mathbb{V}_{\overline{\lambda}} P_{\overline{\lambda}} \lbrace s_{\ld,i}+s_{\ld-\alpha_i,j} \rbrace). $$

Let $ f \in J_{i,j}^-$  and suppose that $ f(m)= \gamma_1 \otimes \gamma_2
\otimes n. $ Then define $ \psi_{-i,-j;\ld} f(m) = \overline{\psi}_{-i,-j;\ld}(\gamma_1 \otimes
\gamma_2) \otimes n. $

\begin{theorem}
\label{catO}
There is a 2-functor $ \Omega_{k,n} \colon \mathcal{KL} \rightarrow \mathcal{P}_{k,n} $ such that for all $i,j\in I$,
\begin{enumerate}
\item $ \Omega_{k,n}(\lambda) =  {}_{\mathbb{Z}} \mathcal{P}_{\overline{\lambda}}^{(k,k)}, $
\item $ \Omega_{k,n}(\mathcal{I}_{\ld)} = \bbI_{\ld}, $
\item $ \Omega_{k,n}(\mathcal{E}_i \mathcal{I}_{\ld}) = \mathbb{E}_i \bbI_{\ld},$
\item $ \Omega_{k,n}(Y_{i;\ld}) = {y}_{i;\ld}, $
\item $ \Omega_{k,n}(\Psi_{i,j;\ld}) = {\psi}_{i,j;\ld},$
\item $ \Omega_{k,n}(\Cup_{i;\ld}) = {\cup}_{i;\ld}, $
\item $ \Omega_{k,n}(\Cap_{i;\ld}) = {\cap}_{i;\ld}, $
\item $ \Omega_{k,n}(\bone_{i;\ld}) = 1_{i;\ld} .$
\end{enumerate}
\end{theorem}

\begin{proof}
This now follows from the computations in \cite[Section 6.2]{kl} for bimodules over the cohomology of flag varieties using the naturality of the isomorphism in proposition ~\ref{fullyfaithful}.
\end{proof}

Finally we show that the category $ \mathcal{P}_{k,n} $ is a categorification of the module $ V_{2 \omega_k}. $
Denote the Grothendieck group of $ \mathcal{P}_{k,n} $ by $ [\mathcal{P}_{k,n}], $ and let $[\mathcal{P}_{k,n}]_{\Q(q)}=\mathbb{C}(q) \otimes_{\mathbb{Z}[q,q^{-1}]}[\mathcal{P}_{k,n}]$.

\begin{prop}
There is an isomorphism of $ \mathcal{U}_q(\mathfrak{sl}_n) $ modules
$[\mathcal{P}_{k,n}]_{\Q(q)} \cong V_{2 \omega_k}. $
\end{prop}

\begin{proof}
Since projective functors map projective-injective modules to projective-injective modules, it follows from
Theorem ~\ref{catO} and \cite{kl}, that $ [\mathcal{P}_{k,n}]_{\Q(q)} $ is a $ \mathcal{U}_q(\mathfrak{sl}_n)$-module.  By construction, it contains a  highest weight vector of weight $ 2 \omega_k $ so it suffices to compute the dimension of its weight spaces.

By \cite[Theorem 4.8]{bk}, the number of projective-injective objects in $\mathcal{O}^{(k,k)}_{\overline{\lambda}}(\mathfrak{gl}_{2k}) $ is equal to the number of column decreasing and row non-decreasing tableau for a diagram with $ k $ rows and $ 2 $ columns with entries from the set
$ \lbrace \underbrace{n, \ldots, n}_{\ld_1}, \ldots, \underbrace{1, \ldots, 1}_{\ld_n} \rbrace. $
Call the set of such tableau $ T. $

Let $ S = \lbrace i \in I^+ | \ld_i =1 \rbrace. $  Denote by $ |S| $ the cardinality of this set.  Consider a Young diagram with $ \frac{|S|}{2} $ rows and $ 2 $ columns.  Let $ T' $ denote the set of tableau on such a column with entries from $ S $ such that the rows and columns are decreasing.  It is well known that the cardinality of the set $ T' $ is the Catalan number $ \frac{\binom{2|S|}{|S|}}{|S|+1}. $
There is a bijection between $ T $ and $ T'. $  For any tableaux $ t' \in T' $ one constructs a tableaux $ t \in T $ by inserting a new box with the entry $ i $ in each column for each $ i \in I^+ $ such that $ \ld_i = 2. $
The inverse is given by box removal.

Finally, the Weyl character formula gives that the dimension of the $ \ld $ weight space of $ V_{2 \omega_k} $ is $ \frac{\binom{2|S|}{|S|}}{|S|+1}. $

\end{proof}


\begin{thebibliography}{99}

\bibitem{bgs} A. Beilinson, V. Ginzburg, W. Soergel, Koszul duality patterns in representation theory, J. Amer. Math. Soc. {\bf 9} (1996), no. 2, 473-527.

\bibitem{b} J. Brundan, Symmetric functions, parabolic category O and the Springer fiber,  Duke Math. J. {\bf 143} (2008), 41-79.

\bibitem{bk} J. Brundan, A. Kleshchev, Schur-Weyl duality for higher levels, Selecta Math. (N.S.) {\bf 14} (2008), no. 1, 1--57.

\bibitem{bk2} J. Brundan, A. Kleshchev, Graded Decomposition Numbers for Cyclotomic Hecke Algebras, arXiv:0901.4450.

\bibitem{bs2} J. Brundan, C. Stroppel,  Highest weight categories arising from Khovanov's diagram algebra II: Koszulity, arXiv:0806.3472.

\bibitem{bs} J. Brundan, C. Stroppel,  Highest weight categories arising from Khovanov's diagram algebra III: category $ \mathcal{O}$, arXiv:0812.1090.

\bibitem{fks} I. Frenkel, Igor, M. Khovanov, Mikhail, C. Stroppel, A categorification of finite-dimensional irreducible representations of quantum $ \mathfrak{sl}_2$ and their tensor products, Selecta Math. (N.S.) {\bf 12} (2006), no. 3-4, 379--431.

\bibitem{hk} S. Huerfano and M. Khovanov, Categorification of some level two representations of $ \mathfrak{sl}_n$, Journal of Knot Theory and Its Ramifications, Vol. 15, No. 6 (2006) 695-713.

\bibitem{k} M. Khovanov, A functor-valued invariant of tangles, Algebr. Geom. Topol. 2 (2002) 665-741

\bibitem{kl1} M. Khovanov and A. Lauda,  A diagrammatic approach to categorification of quantum groups I, arXiv:0803.4121.

\bibitem{kl} M. Khovanov and A. Lauda,  A diagrammatic approach to categorification of quantum groups III, arXiv:0807.3250.

\bibitem{kr} M. Khovanov and L. Rozansky, Matrix factorizations and link homology, Fund. Math. {\bf 199} (2008), no. 1, 1-91.

\bibitem{l} A. Lauda, A categorification of quantum sl(2), math.QA/0803.3652.

\bibitem{m} M. Mackaay,  sl(3)-Foams and the Khovanov-Lauda categorification of quantum sl(k), arXiv:0905.2059.

\bibitem{r} R. Rouquier,  2-Kac-Moody algebras, arXiv:0812.5023.

\bibitem{soe1} W. Soergel, Kategorie $ \mathcal{O} $, perverse Garben und Moduln uber den Koinvarianten zur Weylgruppe, J. Amer. Math. Soc. {\bf 3} (1990), no. 2, 421-445.

\bibitem{soe2} W. Soergel, The combinatorics of Harish-Chandra bimodules, J. Reine Angew. Math. {\bf 429} (1992), no. 2, 49-74.

\bibitem{str1} C. Stroppel, Category $ \mathcal{O}: $ Gradings and Translation functors, Journal of Algebra {\bf 268} (2003), no. 1, 301-326.

\bibitem{str2} C. Stroppel, Categorification of the Temperley-Lieb category, tangles , and cobordisms via projective functors, Duke Math. J. {\bf 126} (2005), no. 3, 547-596.

\bibitem{w} H. Wu, Matrix factorizations and colored MOY graphs, arXiv:0803.2071.

\bibitem{y} Y. Yonezawa, Matrix factorizations and intertwiners of the fundamental representations of quantum group $ \mathcal{U}_q(\mathfrak{sl}_n), $ arXiv:0806.4939.

\end{thebibliography}
\end{document}